\documentclass[11pt,a4paper,reqno]{amsart}
\usepackage{float}
\usepackage{amssymb}
\usepackage{amsmath}
\usepackage{epsfig}
\usepackage[all]{xy}
\usepackage[english]{babel}
\usepackage{latexsym}
\usepackage[T1]{fontenc}
\numberwithin{equation}{section}

\usepackage{ae,aecompl}
\calclayout

\numberwithin{equation}{section}

\renewcommand{\geq}{\geqslant} \renewcommand{\leq}{\leqslant}

\def\TT{{\mathbb{T}}} \def\PP{{\mathbb{P}}} \def\HH{{\mathbb{H}}}
\def\ZZ{{\mathbb{Z}}} \def\QQ{{\mathbb{Q}}} \def\CC{{\mathbb{C}}}
\def\RR{{\mathbb{R}}} \def\EE{{\mathbb{E}}} \def\FF{{\mathbb{F}}}
\def\cO{{\mathcal{O}}}
\def\cF{{\mathcal{F}}}
\def\cN{{\mathcal{N}}}
\def\cP{{\mathcal{P}}}

\renewcommand{\epsilon}{\varepsilon}

\newcommand{\SU}{\mathrm{SU}}

\newcommand{\PU}{\PP\mathrm{U}}
\newcommand{\CP}{\mathbb{CP}}

\newcommand{\sw}{\mbox{}^w}

\renewcommand{\Im}{\mathrm{Im}}

\newcommand{\tU}{\widetilde{U}}
\newcommand{\ovS}{\overline{\Sigma}}
\newcommand{\ovom}{\overline{\omega}}
\newcommand{\ovg}{\overline{g}}
\newcommand{\ovpi}{\overline{\pi}}

\newcommand{\ovF}{\overline{F}}
\newcommand{\ovM}{\overline{M}}
\newcommand{\ovL}{\overline{L}}
\newcommand{\hM}{\widehat{M}}
\newcommand{\chM}{\check{M}}
\newcommand{\chL}{\check{L}}

\newcommand{\ovX}{\overline{X}}

\newcommand{\opi}{\overline{\pi}}
\newcommand{\hS}{\widehat{\Sigma}}
\newcommand{\fa}{\mathfrak a}
\newcommand{\myref}[1]{(\ref{#1})}
\newcommand{\rd}{{\mathrm d}}
\newcommand{\parmu}{\mathrm{par}\mu\;}

\newtheorem{lemma}[subsubsection]{Lemma}
\newtheorem{prop}[subsubsection]{Proposition}
\newtheorem{cor}[subsubsection]{Corollary}

\newtheorem{rmk}[subsubsection]{Remark}
\newtheorem*{rmk*}{Remark}
\newtheorem*{rmks*}{Remarks}
\newtheorem{rmks}[subsubsection]{Remarks}
\newtheorem{theo}[subsubsection]{Theorem}
\newtheorem{theointro}{Theorem}

\newtheorem{conj}[theointro]{Conjecture}
\newtheorem{propintro}[theointro]{Proposition}

\newenvironment{remark*}{\begin{rmk*} --- \normalfont} { \end{rmk*} }
\newenvironment{remarks*}{\begin{rmks*} \begin{enumerate} \normalfont}
{\end{enumerate} \end{rmks*} } 
 
\title{Non-minimal scalar-flat K\"ahler surfaces and parabolic stability}
\date{April 2004}
\author{Yann Rollin}
\thanks{First author partly supported by NSF grant \# DMS-0305130}
\address{Yann Rollin, MIT, 77 Massachusetts Avenue,
Cambridge MA 02139, USA}
\email{rollin@math.mit.edu}
\author{Michael Singer}
\thanks{Second author partly supported by a Leverhulme research
fellowship and an EPSRC small grant}  
\address{Michael Singer, School of Mathematics, King's buildings,
  Edinburgh Scotland} 
\email{msinger@math.mit.edu}
\thanks{}
\begin{document}
{\Huge \sc \bf\maketitle}

\begin{abstract}
A new construction is presented of scalar-flat K\"ahler metrics on non-minimal 
ruled surfaces. The method is based on the resolution of singularities of orbifold
ruled surfaces which are closely related to rank-2 parabolically stable 
holomorphic bundles. This rather general construction is shown also to give new
examples of low genus: in particular, it is shown that 
 $\CP^2$ blown up at $10$ suitably chosen points,  admits a 
scalar-flat K\"ahler metric; this answers a question raised by Claude LeBrun 
in 1986
in connection with the classification of compact self-dual $4$-manifolds.
\end{abstract}

\section{Introduction}

The theory of K\"ahler metrics with constant scalar curvature (CSC) has
seen significant progress in the last ten years: some of the obvious
highlights are:
\begin{itemize}
\item The work of LeBrun and his co-workers \cite{L2,L3, KLP, LS}
which gives many compact examples in complex dimension $2$;

\item The work of Donaldson \cite{D} which shows that for projective
varieties with no non-trivial holomorphic vector fields
there is at most one metric of constant scalar curvature in
any given K\"ahler class; and that if such a metric {\em does} exist,
then the underlying polarized variety must be {\em stable} in a
suitable algebro-geometric sense;

\item The work of X.X. Chen and Tian \cite{CT} which extends Donaldson's
uniqueness result, by different methods, to arbitrary compact 
K\"ahler manifolds (and to {\em
extremal} K\"ahler metrics).
\end{itemize}
The conjecture that the existence of CSC K\"ahler metrics should be
related to the algebro-geometric notion of stability seems to go back
to Yau, but, despite the work of Donaldson, Chen and Tian, a proof of the
`obvious conjecture' (stability $\Rightarrow$ existence) is still
lacking. In the absence of a general theorem of this kind,
special constructions still have an important role to play.

In this paper we give a new construction of compact {\em scalar-flat
K\"ahler} (SFK) surfaces, in other words K\"ahler metrics on compact
complex surfaces having scalar curvature equal to zero. If  $M$ is
such a surface and $[\omega]$ is the K\"ahler class, then
$c_1\cdot[\omega]=0$, for this is just the integral over $M$ of the
scalar curvature. It follows (cf.\ \cite{L1, LS}) that $c_1^2(M) \leq
0$, with equality if and only $c_1(M)=0$ and the metric is
K\"ahler--Einstein. It follows from surface classification 
that if $c_1^2(M)\leq 0$ but
$c_1(M)\not=0$, then $M$ is rational or ruled, though not necessarily
{\em minimal}\footnote{I.e. $M$ could contain divisors that can be
blown down without introducing singularities}. 

The simplest examples are therefore blow-ups of the complex projective
plane $\CP^2$. Since $c_1^2(\CP^2)=9$ and every blow-up reduces
$c_1^2$ by $1$, we see that a 10-point blow-up of $\CP^2$ is the first
such surface that could possibly admit a SFK metric. In \cite{KLP} it
was shown that a $14$-point blow-up of $\CP^2$ does admit a SFK
metric. Our first result is a sharp improvement of this result,
answering a long-standing question \cite{L1} of LeBrun's:

\begin{theointro}
\label{theop2}
The complex projective plane $\CP^2$, blown up at $10$ suitably chosen
points, admits a scalar flat K\"ahler metric. Any further blow-up of
the resulting complex surface admits a scalar-flat K\"ahler metric.
\end{theointro}

We also obtain new constructions of SFK metrics on ruled surfaces with
base an elliptic curve:
\begin{theointro}
\label{theot}
Let $\TT$ be a compact Riemann surface of genus 1.
\begin{enumerate}
\item
Let $L_1$ and $L_2$ be two non-isomorphic holomorphic line bundles of same degree over $\TT$.
Then there is a $2$-point blow-up of $\PP(L_1\oplus L_2)$ which admits a 
scalar-flat K\"ahler metric.
\item
There is a $4$-point blow-up of $\TT\times\CP^1$ which admits a
scalar-flat K\"ahler metric. 
\end{enumerate}
 Any further blow-up of
the resulting complex surfaces admit a scalar-flat K\"ahler metric.
\end{theointro}

In addition to these specific examples, our construction provides some
support for the slogan ``stability $\Rightarrow$ existence''.  For {\em
minimal} ruled surfaces (i.e., no blow-ups) the relation between
stability and existence was noticed in \cite{BB}: a ruled surface of
the form $\PP(E) \to \Sigma$, where $\Sigma$ is a Riemann surface of
genus $\geq 2$ and $E\to
\Sigma$ is a rank-2 holomorphic vector bundle, admits a SFK metric if
and only if $E$ is {\em polystable}. This result depends on the celebrated
theorem of Narasimhan and Seshadri \cite{NS}, which allows to
construct the metric on $\PP(E)$ as a quotient of the
Riemannian product metric $\CP^1\times \HH^2$, where the two factors are
equipped with the standard metrics of constant curvature $+1$ and $-1$
respectively.   

For non-minimal ruled surfaces, the following result was proved by
LeBrun and the second author:

\begin{propintro}[\cite{LS}, Corollary~3.9] Let $M$ be some blow-up of
a compact geometrically ruled surface $\pi:\PP(E)\to \Sigma$. Suppose
that $M$ admits a non-zero, 
periodic holomorphic vector field. Then a K\"ahler class $[\omega]$ on
$M$ with $c_1(M)\cdot [\omega]=0$ contains a representative of zero
scalar curvature if and only if the parabolic bundle $E$ is
quasi-stable.\label{p1.26.11}
\end{propintro}

In this statement, the blow-up and K\"ahler class are encoded by
a parabolic $GL_2(\CC)$-structure on $E$ as follows. If the centres of the
blow-ups are the points $Q_1,\ldots, Q_k$, then we obtain $k$ marked
points $P_j = \pi(Q_j)$ in $\Sigma$ and flags $0 \subset L_j \subset
\pi^{-1}(P_j)$ in the corresponding fibres. The corresponding
parabolic weight $(\beta_j,\gamma_j)$ is not uniquely defined, but
is chosen to satisfy
$$
\gamma_j-\beta_j = \frac{\int_{S_j} \omega}{\int_{F}\omega},\;\;
 \beta_j,\gamma_j \in ]0,1[
$$
where $S_j$ is the exceptional divisor introduced by blowing up
$P_j$.

At the end of \cite{LS} it was conjectured that 
Proposition~\ref{p1.26.11} should continue to 
hold if there is no periodic holomorphic vector field, with
``quasi-stable'' replaced by ``stable''. The methods of this article
do not prove this conjecture; instead, we use the parabolic structure
to encode an iterated blow-up of $\PP(E)$ and define a ``map'' of the
following kind:

\vspace{10pt}

\begin{center}
\begin{tabular}{lcr}
 \fbox{\parbox{2in}{Parabolically stable bundles 
over a punctured hyperbolic Riemann surface $E\rightarrow \widehat \Sigma$}} 
& $\rightarrow$ &\fbox{\parbox{3in}{SFK metrics on a certain blow-up
 $\widehat M$ of $\PP(E)$, encoded by the parabolic structure
}}
\end{tabular}
\end{center}

\vspace{10pt}

\noindent
The reason for the quotation marks is that this ``map'' is only defined
for rational values of the parabolic weights, and will not be smooth
in any obvious sense. {\em We state again that although the blow-up
$\widehat M$ on the right-hand side is encoded by the parabolic
structure, its construction is completely different from the one
involved in Proposition~\ref{p1.26.11}}.

By the theorem of Mehta--Seshadri \cite{MS}, there is a
correspondence
between parabolically stable bundles and representations of the 
fundamental group of
the punctured Riemann surface.  We use these representations in the statement of our
main theorem:
\begin{theointro}
\label{maintheo}
Let $\widehat \Sigma $ be a compact Riemann surface of genus $g$ with a
finite set of marked points $\{P_1,P_2,\cdots, P_k\}$ and  $\rho:
\pi_1(\widehat \Sigma \setminus \{ P_j\}) \rightarrow \SU(2)/\ZZ_2$ be a
homomorphism. Assume in addition that 
\begin{enumerate}
\item  if $l_j$ is the 
homotopy class  of a small loop
around $P_j$,  then $\rho(l_j)$ has finite order $q_j$;
\item $ \displaystyle 2 -2g  - \sum_{j=1}^k (1- \frac 1 {q_j}) < 0$;
\item   $\rho$ defines
  an irreducible representation in the sense that the induced action
  of $\pi_1(\widehat \Sigma \setminus \{ P_j\})$ fixes no point of $\CP^1$.
\end{enumerate}
Then there is non-minimal ruled surface $\hM_\rho \to \hS$ associated
canonically to $\rho$, which admits a
scalar-flat K\"ahler metric.
\end{theointro}

The construction of $\hM_\rho$ will be sketched later in the
Introduction and is given in detail in \S\S2--3. 
\begin{remarks*}
\item We shall see in \S2 that the blow-ups made in the construction
of our SFK ruled surfaces $\hM$ are rather non-generic: they all involve
iterated blow-ups (a sequence of blow-ups where each centre lies on
the exceptional divisor introduced by the previous blow-up). Thus we
get SFK metrics on a rather ``thin'' set in the moduli space of complex
structures on $\hM$.  The general problem of existence of a SFK metric
in a K\"ahler class satisfying $c_1\cdot [\omega]=0$ remains
mysterious.

\item It is not clear whether Theorem~B is sharp: for example, do
there exist SFK metrics on a 1-point blow-up of $\PP(L_1\oplus L_2)\to
\TT$, if $L_1$ and $L_2$ are not isomorphic? Do there exist such
metrics on a 3-point blow-up of $\TT \times \CP^1$?  In this direction,
we remark that \cite[Prop.~3.1]{LS} shows that there do not exist SFK
metrics in the first case if the centre of the blow-up lies on $L_1$
or $L_2$ and in the second case if 2 or fewer points are blown
up. (The obstructions come from the non-trivial holomorphic vector
field on these spaces.)
\end{remarks*}

In the light of the recent work of Donaldson and Chen--Tian, the
following seems a reasonable

\begin{conj}
Let $E\to \Sigma$ be a parabolic holomorphic bundle of rank $2$ (with
rational weights)  over a
Riemann surface, such that the ruled surface $\PP(E)$ has  no
non-trivial holomorphic vector field. If the corresponding iterated
blow-up $\widehat M$  
of $\PP(E)$ (cf. Section~\ref{secitbup} for a precise definition) admits 
a scalar-flat K\"ahler metric, then $E$ must be parabolically stable.  
\end{conj}

It is ironic that the above conjecture runs in the ``easier'' direction
(existence $\Rightarrow$ stability) and yet we are unable to prove
it; the new numerical criterion for stability due to Ross and Thomas
\cite{RT} should be useful here, but so far we have not succeeded in
applying it.

Another possible approach might be to apply the Tian--Viaclovsky
compactness theorem \cite{TV} which shows that under certain
conditions, a sequence of SFK metrics can only degenerate to a SFK
orbifold metric.  (Our gluing theorem, Theorem~\ref{theoglue}, gives
explicit examples of this degeneration process.)

\subsection*{Outline}

This work began from the observation that in the Burns--de Bartolomeis
construction, the smooth base $\Sigma$ can be replaced by  an orbifold
Riemann surface $\ovS$. More precisely, with the notation of
Theorem~\ref{maintheo}, $\ovS$ is the smooth Riemann surface $\Sigma$
with a finite set of marked points $P_j$ and corresponding integer
weights $q_j$.  By a theorem of Troyanov \cite{Tr}, if the orbifold
Euler characteristic is negative (this is condition (ii) of
Theorem~\ref{maintheo}), then $\ovS$ carries a (K\"ahler) orbifold
metric $\overline{g}$ of constant curvature $-1$. In particular, 
$\overline{g}$ is smooth on $\ovS \setminus\{P_j\}$ and has a conical
singularity at 
$P_j$, with cone angle $2\pi/q_j$. The riemannian product $\ovS \times
\CP^1$ is obviously scalar-flat K\"ahler, with non-isolated orbifold
singularities around the fibres $F_j = \{P_j\}\times \CP^1$.

In order to replace these by isolated singularities, we twist by a
representation $\rho$ of the orbifold fundamental group of $\ovS$,
just as was done in the smooth case by Burns and de Bartolomeis.  This
$\rho$ must be as in condition (i) of Theorem~\ref{maintheo}: using it
gives an orbifold SFK metric on an orbifold ruled surface, $\ovM$,
say. It turns out that $\ovM$ has precisely two isolated cyclic
singularities in each fibre $F_j$. Denote by $\hM$ the minimal
resolution of singularities of $\ovM$.

In \S\S4--5, it will be shown that $\hM$ has a SFK metric, by an
analytical gluing theorem. On the other hand, we shall see in \S\S2--3
that $\hM$ is also a multiple blow-up of a smooth minimal ruled
surface $\chM= \PP(E)\to \hS$, say. In fact, $\chM$ and $\ovM$ can
be viewed as two different compactifications (one smooth, the other an
orbifold) of a non-compact ruled surface $M^* \to
\Sigma\setminus\{P_j\}$, each canonically associated to the
representation $\rho$. Such representations are related by the
Mehta--Seshadri Theorem \cite{MS} to parabolic stability of the
underlying holomorphic vector bundle.  In \S2 we shall start from this
point, defining a notion of parabolically stable ruled surface $\chM$
and the corresponding multiple blow-up $\hM$. In \S3, we shall compare
this with the orbifold $\ovM$.

In \S2, we shall restate Theorem~\ref{maintheo} in the language of stable 
ruled surfaces and show how Theorems A and B follow.  The advantage of working with parabolically stable bundles is that it is often quite easy to verify stability, whereas it can be rather difficult to find explicit representations of the fundamental group of a punctured
Riemann surface.   

%
%


\subsection*{Acknowledgments}

We thank Olivier Biquard for useful conversations and for
pointing out a mistake in the proof of the gluing theorem for
scalar-flat K\"ahler metrics in 
\cite{KS}. This problem in now fixed in section \ref{secgluing} where a
self-contained exposition of the result can be found.  We also thank
Claude LeBrun and Rafe Mazzeo for encouragement and several useful
discussions. 
Finally, we thank the anonymous referee for 
several suggestions which improved the original manuscript.
This work was carried out while the
second author was visiting the Mathematics Department at MIT; he
thanks MIT for its hospitality and financial support
during this visit.

\section{Parabolically stable ruled surfaces}
\label{secparab}
A {\em geometrically ruled surface} $\check M$ is by definition a minimal complex
surface obtained as $\check M= \PP(E)$, where $E\rightarrow \widehat \Sigma
$ is a holomorphic vector bundle of rank $2$ over a Riemann surface 
$\widehat \Sigma$.
The induced map $\pi:\check M \rightarrow \widehat \Sigma$ is called the
{\em ruling}.

A parabolic structure on $\check M$ consists of the following data:
\begin{itemize}
\item A finite set of distinct points $P_1,P_2,\cdots,P_n$ in
  $\widehat \Sigma$;
\item for each $j$, a choice of point  $Q_j \in F_j = \pi^{-1}(P_j)$;

\item for each $j$, a choice of  weight $\alpha_j \in   ]0,1[\cap\QQ$.
\end{itemize}
A geometrically ruled surface with a parabolic structure will be called a
{\em parabolic ruled surface}.

If $S\subset\check M$ is a holomorphic section of $\pi$, we define its
slope
$$\mu(S) = S^2 +\sum_{Q_j\not\in S}\alpha_j - \sum_{Q_j\in S}\alpha_j;
$$
we say that a parabolic ruled surface is \emph{stable} if for
every holomorphic section $S$, we have $\mu(S) >0$.

If we return to the vector bundle $E$, 
then $Q_j$ defines a line $L_j$ in the fibre of $E$ over $P_j$. For
each $j$, select $0\leq  \beta_j < \gamma_j < 1$ with
$\alpha_j = \gamma_j- \beta_j$. In this way $E$ is endowed with (a
family of) parabolic structures. Our notion of stability of a
parabolic ruled surface corresponds with the Mehta--Seshadri notion of
parabolic stability for $E$. Indeed, holomorphic cross-sections
$S\subset\chM$ correspond exactly to holomorphic sub-bundles $L\subset
E$. 
We know that 
$H^2(\PP(E),\ZZ)$ is generated by the class of a fibre $F$ and
 $H = c_1(\cO(1))$ (the fibrewise hyperplane section bundle) on each fibre (cf. for instance
\cite{Beau}). They verify 
$$ H^2=\deg(E),\quad F^2=0,\quad F\cdot H=1.
$$
Moreover, we have 
$$ S = H - \deg(L) F
$$
hence
$$ S^2 = \deg(E)- 2\deg(L).
$$
It follows that $\mu(S)$ is equal to twice the difference of the
parabolic slopes of $E$ and $L$ in the sense of
Mehta--Seshadri~\cite{MS}, so 
$\check M$ is stable if and only if $E$ is parabolically stable.

\subsection{Iterated blow-up of a parabolic ruled surface}
\label{secitbup}

Let $\chM$ be a parabolic ruled surface. We shall now define a
multiple blow-up $\Phi:\hM \to \chM$ which is canonically
determined by the parabolic structure of $\chM$.

In order to simplify the notation, suppose that the parabolic
structure on $\check M$ is reduced to a single point $P\in \widehat
\Sigma$; let $Q$ be the corresponding point in $F = \pi^{-1}(P)$ and let
$\alpha = \frac pq$ be the weight, where $p$ and $q$ are two coprime
integers, $0<p<q$. Denote the Hirzebruch--Jung continued fraction
expansion of $\alpha$ by
\begin{equation}\label{e1.844}
\frac pq =  \cfrac{1}{e_1-\cfrac{1}{e_2-\cdots\cfrac{1}{e_k}}};
\end{equation}
define also
\begin{equation}\label{e20.844}
\frac {q-p}{q} = \cfrac{1}{e'_1-\cfrac{1}{e'_2-\cdots
\cfrac{1}{e'_l}}}.
\end{equation}
These expansions are unique if, as we shall assume, the $e_j$ and $e'_j$ 
are all $\geq 2$.

\begin{prop} There exists a unique iterated blow-up $\Phi:\hM\to\chM$
with $\Phi^{-1}(F)$ equal to the following chain of curves:
\begin{equation}
\xymatrix{
{}\ar@{-}[r]^{-e_1} & *+[o][F-]{}
\ar@{-}[r]^{ -e_{2}} &  *+[o][F-]{}
\ar@{--}[r] &  *+[o][F-]{}
\ar@{-}[r]^{-e_{k-1}} &  *+[o][F-]{}
\ar@{-}[r]^{-e_k} &  *+[o][F-]{}
\ar@{-}[r]^{-1} &  *+[o][F-]{}
\ar@{-}[r]^{-e'_{l}} &  *+[o][F-]{}
\ar@{-}[r]^{-e'_{l-1}} &  *+[o][F-]{}
\ar@{--}[r] &  *+[o][F-]{}
\ar@{-}[r]^{-e'_{2}} &  *+[o][F-]{}
\ar@{-}[r]^{ -e'_{1}} & 
}.
\end{equation}

\noindent
Here the edges represent rational curves, the number above each edge  
is the self-intersection of the curve, 
the hollow dots represent transverse intersections with intersection number
$+1$, and the
curve of self-intersection $-e_1$ is the proper transform of the
exceptional divisor of the first blow-up. 
\label{p1.16.4.4}\end{prop}

\begin{proof} We give an iterative construction. The first step is to
blow up $Q$, to get a diagram of the form
\begin{equation}\label{e2.844}
\xymatrix{
{}\ar@{-}[rr]^{-1} && *+[o][F-]{}
\ar@{-}[rr]^{-1} &&{}
}
\end{equation}
By blowing up the intersection point of these two curves we get the
diagram
\begin{equation}\label{e3.844}
\xymatrix{
{}\ar@{-}[rr]^{-2} && *+[o][F-]{} 
\ar@{-}[rr]^{-1} && *+[o][F]{}
\ar@{-}[rr]^{-2} &&{} 
}
\end{equation}
in which we see two $(-2)$-curves separated by a $-1$ curve. Suppose
by induction that we have used a sequence of blow-ups so that the
following chain of curves sits over $F$:
\begin{equation}\label{e4.844}
\xymatrix{
{}\ar@{-}[r]^{-a_1} & *+[o][F-]{}
\ar@{-}[r]^{ -a_{2}} &  *+[o][F-]{}
\ar@{--}[r] &  *+[o][F-]{}
\ar@{-}[r]^{-a_{j-1}} &  *+[o][F-]{}
\ar@{-}[r]^{-a_j} &  *+[o][F-]_{A}{}
\ar@{-}[r]^{-1} &  *+[o][F-]_{B}{}
\ar@{-}[r]^{-b_{r}} &  *+[o][F-]{}
\ar@{-}[r]^{-b_{r-1}} &  *+[o][F-]{}
\ar@{--}[r] &  *+[o][F-]{}
\ar@{-}[r]^{-b_{2}} &  *+[o][F-]{}
\ar@{-}[r]^{ -b_{1}} & 
}
\end{equation}

\noindent
for integers $a_i$ and $b_i$ $\geq 2$. Then we can increase $a_j$ by
one unit by blowing up the point marked $A$ and we can introduce a new
curve of self-intersection $-2$ by blowing up $B$:
\begin{equation}\label{e5.844}
\xymatrix{
{}\ar@{-}[r]^{-a_1} & *+[o][F-]{}
\ar@{-}[r]^{ -a_{2}} &  *+[o][F-]{}
\ar@{--}[r] &  *+[o][F-]{}
\ar@{-}[r]^{-a_{j-1}} &  *+[o][F-]{}
\ar@{-}[r]^{-a_j-1} &  *+[o][F-]{}
\ar@{-}[r]^{-1} &  *+[o][F-]{}
\ar@{-}[r]^{-2} &  *+[o][F-]{}
\ar@{-}[r]^{-b_{r}} &  *+[o][F-]{}
\ar@{-}[r]^{-b_{r-1}} &  *+[o][F-]{}
\ar@{--}[r] &  *+[o][F-]{}
\ar@{-}[r]^{-b_{2}} &  *+[o][F-]{}
\ar@{-}[r]^{ -b_{1}} & 
}
\end{equation}

or

\begin{equation}\label{e6.844}
\xymatrix{
{}\ar@{-}[r]^{-a_1} & *+[o][F-]{}
\ar@{-}[r]^{ -a_{2}} &  *+[o][F-]{}
\ar@{--}[r] &  *+[o][F-]{}
\ar@{-}[r]^{-a_{j-1}} &  *+[o][F-]{}
\ar@{-}[r]^{-a_j} &  *+[o][F-]{}
\ar@{-}[r]^{-2} &  *+[o][F-]{}
\ar@{-}[r]^{-1} &  *+[o][F-]{}
\ar@{-}[r]^{-b_{r}-1} &  *+[o][F-]{}
\ar@{-}[r]^{-b_{r-1}} &  *+[o][F-]{}
\ar@{--}[r] &  *+[o][F-]{}
\ar@{-}[r]^{-b_{2}} &  *+[o][F-]{}
\ar@{-}[r]^{ -b_{1}} & 
}
\end{equation}
It is now clear that a suitable sequence of $\sum(e_j -1)$ blow-ups,
starting from the configuration \eqref{e2.844} yields the required
configuration of curves. All that remains is the proof that the integers
$e'_j$ in \eqref{e2.844} do satisfy \eqref{e20.844}. For this we use the
fact if
\begin{equation}\label{e15.844}
\frac{\lambda}{\mu} = 
\cfrac{1}{a_1-\cfrac{1}{a_2- \cdots \frac{1}{a_j}}}
\end{equation}
then the reversed continued fraction
\begin{equation}\label{e16.844}
\cfrac{1}{a_j-\cfrac{1}{a_{j-1}- \cdots \frac{1}{a_1}}}
= \frac{\lambda'}{\mu}
\end{equation}
where $\lambda \lambda' \equiv 1 \pmod \mu$, $0< \lambda' <
\mu$. Suppose by induction that in \eqref{e4.844} we have
\begin{equation}\label{e7.844}
\cfrac{1}{a_1-\cfrac{1}{a_2- \cdots \frac{1}{a_j}}} +
\cfrac{1}{b_1-\cfrac{1}{b_2- \cdots \frac{1}{b_r}}} = 1
\end{equation}
Since $(\mu - \lambda')(\mu - \lambda)$ is congruent to 1 mod $\mu$,
\eqref{e7.844} implies that
\begin{equation}\label{e8.844}
\cfrac{1}{a_j-\cfrac{1}{a_{j-1}- \cdots \frac{1}{a_1}}} +
\cfrac{1}{b_r-\cfrac{1}{b_{r-1}- \cdots \frac{1}{b_1}}} = 1
\end{equation}
Now consider how these fractions change under the blow-up. We have
that $a_j$ is replaced by $a_j+1$, and a new $b_{r+1}$ is introduced,
equal to $2$. So the
inductive step is completed by the trivial identity
$$
\frac{1}{1 + \mu/\lambda'} + \frac{1}{2- (1-\lambda'/\mu)} = 1.
$$
Since the induction clearly starts (consider the diagram
\eqref{e3.844}), the proof is now complete.
\end{proof}

Our main theorem may now be stated as follows:
\begin{theo}
\label{maintheoparab}
Let $\check M\rightarrow \widehat \Sigma$ be a parabolically stable
ruled surface. Suppose that 
\begin{equation}
\label{eqhyp}
\chi(\widehat \Sigma) - \sum_{j=1}^k(1 - \frac 1{q_j}) <0,
\end{equation}
where $\alpha_j = p_j/q_j$, with $p_j$ and $q_j$ coprime.
Then the iterated blow-up $\widehat M$ deduced from the parabolic
structure carries a SFK metric.  Furthermore, any blow-up of $\hM$ 
also carries a SFK metric.
\end{theo}

The relationship between this statement and Theorem~\ref{maintheo}
corresponds to the relationship between representations and
parabolically stable bundles, given by the 
theorem of Mehta and Seshadri \cite{MS}. According to this theorem,
given a parabolically stable $\chM$, there is a corresponding
irreducible representation
\begin{equation}\label{e1p.16.4.4}
\rho' :\pi_1(\hS \setminus \{P_j\}) \to U_2
\end{equation}
with
\begin{equation}\label{e2p.16.4.4}
\rho'(l_j) \mbox{ conjugate to }\pm \mbox{diag}(e^{2 \pi i\beta_j},
e^{2\pi i\gamma_j})
\end{equation}
where $l_j$ is the homotopy class of a small loop around
$P_j$ and $\gamma_j-\beta_j= \alpha_j$ as before. This representation
defines a parabolically stable bundle $E$ with $\PP(E) = \chM$. We
note that $E$ is determined only up to tensoring with a line-bundle;
we make use of this freedom (and the freedom of choice in $\beta_j$,
$\gamma_j$) to ensure that the parabolic degree of  $E$ is zero, as is
necessary for the existence of $\rho'$.

If $A\in U_2$, write
$A = z(A)A_0$, where $z(A)$ is a scalar multiple of the identity and
$A_0\in SU_2/\ZZ_2$. Letting $\rho(\gamma) = \rho'(\gamma)_0$, we
obtain an irreducible representation
\begin{equation}\label{e1.16.4.4}
\rho :\pi_1(\hS \setminus \{P_j\}) \to SU_2/\ZZ_2
\end{equation}
with
\begin{equation}\label{e2.16.4.4}
\rho(l_j) \mbox{ conjugate to }\pm \mbox{diag}(e^{\pi i\alpha_j},
e^{-\pi i\alpha_j})
\end{equation}
as required for Theorem~\ref{maintheo}.

Conversely, given a representation with
properties \eqref{e1.16.4.4} and \eqref{e2.16.4.4}, 
we clearly get a ruled surface $M^* \to \hS
\setminus\{P_j\}$. However, $M^*$ can now be compactified to give
$\chM$ in a standard way: let $P$ stand for one of the punctures, and
identify a neighbourhood of $P$ in $\hS$ with the unit disc
$\Delta$. Write $L^*= \pi^{-1}(\Delta\setminus\{0\})$ and note that
\begin{equation}\label{e30.844}
L^* = \HH^2\times \CP^1/ \ZZ
\end{equation}
where $\HH^2$ is the hyperbolic plane $\{\Im(\xi) > 0\}$ and the
$\ZZ$-action is generated by
\begin{equation}\label{e31.8.44}
(\xi,[w_0:w_1]) \mapsto (\xi + 2\pi, [e^{\pi i \alpha}w_0,
e^{-\pi i \alpha}w_1]).
\end{equation}
Let $\chL = \Delta \times \CP^1$ be covered by two standard
holomorphic coordinate charts, $(x_1,y_1)\in \CC^2$, $(x_2,y_2)\in
\CC^2$ glued together by $x_2= x_1$, $y_2 = y_1^{-1}$. Then the
natural embedding of $L^*$ in $\chL$ is given by the following map
\begin{equation}\label{e31.844}
x_1 = x_2 =e^{i\xi},\; y_1 = e^{-i\xi \alpha}w_0/w_1,\; 
y_2 = e^{i\xi \alpha}w_1/w_0.
\end{equation}
Notice that the pre-image of $Q$, which corresponds to $(x_1,y_1) =
(0,0)$ is given by $\xi = +i\infty$, $w_0=0$, but every other point
corresponding $\xi= i\infty$ is mapped to $(x_2,y_2) = (0,0)$.  

In general, we obtain
$\chM$ as a smooth compactification of $M^*$ by following this recipe 
at each puncture.

\vspace{12pt}

Let us show how Theorem~\ref{theop2} and Theorem~\ref{theot} follow
from Theorem~\ref{maintheoparab}.

\begin{cor}
There is a $9$-point iterated blow-up of $\CP^1\times\CP^1$ which
admits a SFK 
metric. As a consequence $\CP^2$  has a $10$-point iterated blow-up 
which admits a scalar-flat K\"ahler metric.
\end{cor}
\begin{proof} Let $\chM = \CP^1\times \CP^1$, and let
$\pi:\chM \rightarrow \CP^1$ denote
projection on the first factor. Pick any
$4$ points $P_1,P_2,P_3,P_4$ in $\CP^1$, with weights
$\alpha_1=\alpha_2= \alpha_3 = 1/2$ and $\alpha_4=1/3$, so that
(\ref{eqhyp}) is satisfied. Pick $Q_j \in \pi^{-1}(P_j)$. To check
when this parabolic structure is stable, note first that any section
$S$ is in this case the graph of a meromorphic function $f$, and 
$S^2=2\deg f$. Provided that no two of the $Q_j$ lie on the graph of a
function of degree $0$,
we have a parabolically stable ruled surface. In particular, the set
of stable configurations of the $Q_j$ is a Zariski open set in the set
of all such configurations.

The multiple blow-up $\hM$ of $\chM$ involves a total of 9 blow-ups, 2
each at $Q_1$, $Q_2$ and $Q_3$ and $3$ at $Q_4$. Hence $\hM$ is indeed
a $9$-point blow-up of $\CP^1\times\CP^1$. However, if $n\geq 1$, then
any $n$-point blow-up of $\CP^1\times\CP^1$ is isomorphic to an
$(n+1)$-point blow-up of $\CP^2$.
\end{proof}

Now we turn to ruled surfaces with base an elliptic curve.

\begin{cor}
\label{cortt}
Let $\TT$ be an elliptic curve (Riemann surface of genus 1) and let
$L_1$ and $L_2$ be two non-isomorphic line-bundles of degree $0$ over
$\TT$. Then, any double blow-up 
$\widehat M$ of $\chM=\PP(L_1\oplus L_2)$ at a point which is not on $\PP(L_1)$
or $\PP(L_2)$  admits a SFK metric. Any further blow-up of $\widehat 
M$ also admits a SFK metric.
\end{cor}
\begin{proof} Endow $\chM$ with a parabolic structure by picking an
arbitrary point $P\in \TT$ and a point $Q \in \pi^{-1}(P)$ so that $Q$
does not lie on $S_1=\PP(L_1 \oplus 0)$ or $S_2=\PP(0\oplus
L_2)$. Set the weight of $Q$ equal to $1/2$. Then
condition~(\ref{eqhyp}) is satisfied and it remains to check that this
parabolic structure is stable.

To see this, note first that every section $S$ of $\pi$ satisfies
$S^2\geq 0$; this follows because at the level
of cohomology, $[S] = [S_1] + r[F]$ for some integer $r$. Then $S^2=0$
if $r=0$; otherwise $S\not= S_1$, so $r =S\cdot S_1 \geq 0$ Hence $S^2
= 2r \geq 0$. 

 Therefore, $\mu(S)>0$ unless 
$S^2=0$ and $Q\in S$. Suppose now that  $S^2=0$, $S\neq S_1$ and  $S\neq S_2$ .
It follows that $S\cdot S_1= S\cdot S_2=0$, and $S$ meets neither
$S_1$ nor $S_2$. Such an $S$ defines an isomorphism
$L_1 \simeq L_2$ contradicting the hypothesis of the corollary.

The corresponding 2-point blow-up $\hM$ of $\chM$ carries a
scalar-flat K\"ahler metric.
\end{proof}

In the next corollary, we recover a result of Kim, LeBrun and
Pontecorvo~\cite{KLP} by a completely different method.
\begin{cor}
\label{corklp}
Let $\TT$ be an elliptic Riemann surface. There is a $6$-point blow-up
of $\TT\times \CP^1$ which admits a SFK metric.
\end{cor}
\begin{proof} Let $\chM = \TT\times \CP^1$, let $\pi$ be the
projection on the first factor. Endow $\chM$ with a parabolic
structure by choosing any $3$ points $P_j$ in $\TT$, points $Q_j \in
\pi^{-1}(P_j)$ each of weight $1/2$. We note that the condition
\eqref{eqhyp} is satisfied as before: it remains to analyze the
stability of this parabolic structure. 

Any holomorphic section $S$ of $\pi$ corresponds 
to the graph of a meromorphic function  $s$ on $\TT$. 
Let $H$ be the section
corresponding to the constant function $0$. Then $S\cdot F=1$, and
$S\cdot H =\sum a_n$, where the $a_n\geq 0$ are the multiplicities of
the zeroes of $s$;
hence
$$ S= H + (\sum a_n) F. 
$$  

If $S^2=2 \sum a_n =0$, then $a_n=0$ and $s$ must be  constant. If we
choose the points $Q_j$ such that no two of them lie 
on the graph of a constant meromorphic function, we must have
$\mu(S) \geq 0 + 2/2> 0$.  If $S^2>0$, we have $S^2\geq 2$ for $S^2$
is even and then $\mu(S)\geq 2-3/2>0$ whatever the positions of the
$Q_j$.

Applying Theorem~\ref{maintheoparab}, there exists a SFK metric on
$\hM$, which is in this case a $6$-point blow-up of $\chM$; more
precisely $\chM$ is obtained by performing a double blow-up at each of
$Q_1$, $Q_2$ and $Q_3$.
\end{proof}

We now improve the last result.
\begin{cor}
Let $\TT$ be an elliptic curve.
There is a $4$-point blow-up of $\TT\times \CP^1$ which admits a scalar-flat 
K\"ahler metric.
\end{cor}
\begin{proof} 
Let $P_1$ and $P_2$ be two points in $\TT$.  Let $S$ be a
constant section  which is not the zero section $S_0$ or the section at
infinity $S_\infty$, and $M'$ be the blow-up of $\TT\times\CP^1$ at   
$\{P_1\}\times\{0\}$ and $\{P_2\}\times\{\infty\}$. In $M'$, the proper 
transform of $\{P_j\}\times \CP^1$ is a $(-1)$-curve; we blow each of these curves down, getting a new minimal ruled surface $\chM$, giving the following diagram
$$
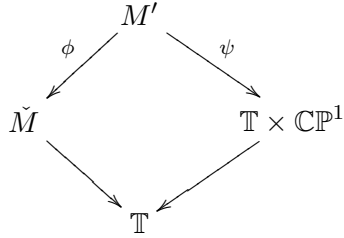
\xymatrix{&          M'\ar[dl]_\phi \ar[dr]^{\psi}   &  \\
 \chM \ar[dr]& & \TT\times\CP^1 \ar[dl] \\
 &\TT & }
$$
Here $\phi$ is a blow-up map with centres at $C_1$ and $C_2$, say, and 
$\phi^{-1}(C_j)$ is the $\psi$-proper transform of $\{P_j\}\times \CP^1$.

We shall show that $\chM = \PP(L_1\oplus L_2)$ satisfies the
conditions of Corollary~\ref{cortt}. Then the double blow-up 
$\Phi:\hM \to \chM$ at some point $Q$ is scalar-flat K\"ahler, and so is
any further blow-up of $\hM$.  
In particular, the blow-up $\widehat{\phi}:\hM'\to \hM$ of $\hM$ 
at $\Phi^{-1}(C_1)$  and $\Phi^{-1}(C_2)$ 
is scalar-flat K\"ahler. Taking $Q$ to be different from 
$C_1$ and $C_2$, we see that $\hM'$ is also the double blow-up at
$\phi^{-1}(Q)$ of $M'$, 
as in the following diagram:
$$
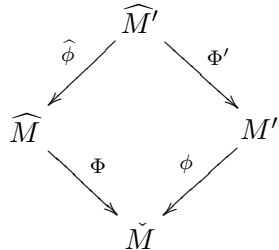
\xymatrix{&          \hM'\ar[dl]_{\widehat{\phi}} \ar[dr]^{\Phi'}   &  \\
 \hM \ar[dr]^{\Phi}& & M' \ar[dl]_\phi \\
 &\chM & }
$$
In particular, $\hM'$ is a $4$-point blow-up of $\TT\times \CP^1$.

Let us check that $\chM$ satisfies the hypotheses of
Corollary~\ref{cortt}. Suppose not. 
Let $S', S'_0, S'_\infty$ be the images by $\phi$ in $\chM$ of the 
proper transforms of the corresponding 
curves $S$, $S_1$ and $S_\infty$ in $\TT\times \CP^1$. We have
$$ (S_0')^2= (S'_\infty)^2= S_0'\cdot S_\infty' =0,\quad S_0'\cdot S'
= S_\infty'\cdot S'=1,\quad (S')^2 = 2.
$$
This shows that $M =\PP(L_0\oplus L_\infty)$, where $L_0$ and
$L_\infty$ are two line bundles of degree $0$ over $\TT$.

Suppose that $M\simeq \TT\times \CP^1$. Then $S_0'$ and $S_\infty'$
must be two distinct constant sections of $\TT\times \CP^1\to\TT$ as
seen in the proof of Corollary~\ref{corklp}. Up to an ismorphism of
$\CP^1$, we may assume that they are the $0$ section and the section
at infinity of $\TT\times\CP^1$. Now $S'$ is the graph of a
meromorphic function $s'$ on $\TT$. Let $z_j$ be the zeroes of $s'$
and $w_j$ its poles. The divisor of $s'$ is $$ (s')= \sum a_j z_j -
\sum b_j w_j, $$ where $a_j, b_j\geq0$ are the multiplicities.  Then,
$S'\cdot S'_0= \sum a_j=1$ hence there is a unique pole with
multiplicity $1$. Therefore $$(s')= z-w, $$ which is impossible by
Abel's theorem.
\end{proof}

\begin{remark*}
It is a general principle in algebraic geometry that ``stability is an
open condition''. We have seen in these examples that the set of stable
parabolic structures is Zariski dense in the set of all such
structures. This will be true in general: if stable parabolic
structures exist on a given ruled surface, then they will form a
Zariski-open subset in the set of all such structures.  Furthermore,
stability will be preserved under perturbation of the parabolic
weights, keeping the $Q_j$ fixed. However the iterated blow-up $\hM$
will behave rather wildly under such perturbations: for example,
$b_2(\hM)$ will not remain constant.
\end{remark*}

\section{Scalar-flat K\"ahler orbifolds}

The goal of this section is to introduce an orbifold $\ovM$ associated
to a parabolic ruled surface $\chM$. We shall see that $\chM$ and
$\ovM$ are bimeromorphic to each other, and indeed both are
compactifications of the punctured ruled surface $M^*$. We shall also
show that $\ovM$ carries an orbifold SFK metric if $\chM$ is stable,
and that $\ovM$ carries no non-trivial holomorphic vector fields.

\subsection{Generalities}
A complex {\em orbifold} $X$ of complex dimension $n$ may be defined
as a complex variety having only quotient singularities. More
explicitly, for every point $P\in X$, one requires that there is a
finite group $G= G_P$ (called the {\em local isotropy group}) and a
{\em local uniformizing chart} 
\begin{equation}\label{e1.10.12.03}
\tU \to \tU/G \stackrel{\phi}{\longrightarrow} U.
\end{equation}
Here, $\tU$ is a neighbourhood of 
$0$ in $\CC^n$, with a given biholomorphic action of $G$,
\begin{equation}\label{e1.23.12.03}
G\times \tU \to \tU,
\end{equation}
 $U$ is
a neighbourhood of $P$ in $X$ and $\phi$ is a homeomorphism with
$\phi(0) = P$. It is a fact that one can always choose $\phi$ so as to linearize the action \eqref{e1.23.12.03}; that is, one can assume
$G\subset GL_n(\CC)$.

The standard notions of differential geometry extend to orbifolds by
working $G$-equivariantly in a local uniformizing chart. For example
a smooth orbifold Riemannian metric $g$ on $X$ is defined as usual
away from the singular points, and is given by a $G_P$-invariant smooth
metric on $\tU$ near the singular point $P$.

\subsection{Orbifold Riemann surfaces}
A compact orbifold Riemann surface $\ovS$ can be identified with a smooth
compact Riemann surface $\hS$ together with a finite set of marked points $P_j$, each with a given weight $q_j\in \ZZ_{\geq 2}$.  It is important to note that a smooth orbifold metric on $\ovS$ is not the same as a smooth metric on $\hS$: a smooth orbifold metric is smooth on $\Sigma = \ovS \setminus\{P_j\}$ but with respect to such a metric, 
the length of a small circle of radius $r$
centred at $P_j$, will be approximately $2 \pi r/q_j$.

The Euler characteristic of $\ovS$ is defined as follows
\begin{equation}\label{e2.10.12.03}
\chi(\ovS) = \chi(\hS) - \sum \left(1 - \frac{1}{q_j}\right),
\end{equation}
where $\chi(\hS)$ is the Euler characteristic of the underlying smooth surface
$\hS$.

Just as for smooth Riemann surfaces with negative Euler characteristic, we have the following result of Troyanov \cite{Tr} (see also \cite{McO}).

\begin{theo} The orbifold Riemann surface $\ovS$ admits an orbifold metric of constant curvature $-1$ compatible with the given complex structure if and only if $\chi(\ovS) <0$.\label{theotr}
\end{theo}

In view of this theorem we shall refer to such orbifold Riemann surfaces as {\em hyperbolic}. In the previous section we saw several examples of orbifold hyperbolic Riemann surfaces.

Next we come to the fundamental group of an orbifold Riemann
surface. Recall first the description of the fundamental group of the 
punctured Riemann surface $\Sigma$:
\begin{equation}\label{e3.10.12.03}
\pi_1(\Sigma)= \langle a_1,b_1,\ldots,a_g,b_g,l_1,\ldots l_k:
[a_1,b_1][a_2,b_2]\ldots[a_g,b_g]l_1\ldots l_k = 1\rangle
\end{equation}
Here the $a_j$ and $b_j$ are standard generators of $\pi_1(\hS)$
 and  $l_j$ is
(the homotopy class of) a small loop around $P_j$. The orbifold
fundamental group is defined by imposing the additional conditions
\begin{equation}\label{e4.10.12.03}
\pi^{orb}_1(\ovS) = \langle 
 a_1,b_1,\ldots,a_g,b_g,l_1,\ldots l_k:
[a_1,b_1][a_2,b_2]\ldots[a_g,b_g]l_1\ldots l_k =  l_1^{q_1}
= \ldots = l_k^{q_k} = 1
\rangle
\end{equation}
From Theorem~\ref{theotr} it follows that for any 
orbifold hyperbolic Riemann surface $\ovS$, we have $\ovS =  \HH^2/\Gamma$, where 
$\Gamma$ is the image of the uniformizing representation of $\pi_1^{orb}(\ovS)$ in $SL_2(\RR) = {\rm Isom}(\HH^2)$. The only difference from the smooth case is that $\Gamma$ will not  act freely on $\HH^2$.

\subsection{Orbifold ruled surfaces}

Let $\ovS$ be a hyperbolic orbifold Riemann surface as in the last section, and let 
$\rho: \pi_1^{orb}(\ovS) \to PSL_2(\CC)$ be a representation. We suppose that $\rho$ is faithful on the loops $l_j$, so that $\rho(l_j)$ has order precisely $q_j$. We can form the quotient
\begin{equation}\label{e1.154}
\ovM = \HH^2\times \CP^1/\pi_1^{orb}(\ovS)
\end{equation}
by letting $\pi_1^{orb}(\ovS)$ act by the uniformizing representation
on the upper half-space $\HH^2$ 
and by $\rho$ on $\CP^1$. It is clear that $\ovM$ is a
complex orbifold ruled surface equipped with a ruling 
$$
\ovpi: \ovM \to \ovS
$$
and singularities only in the fibres $\ovF_j=\ovpi^{-1}(P_j)$.

In order to analyze these singularities, choose a complex disc
$\Delta$ with centre at one of the $P_j$, let $\ovL$ denote
$\ovpi^{-1}(\Delta)$ and let $L^* = \ovL \setminus \ovpi^{-1}(0)$.
Complex analytically, we can identify $L^*$ with a quotient as before
(see \eqref{e30.844} and \eqref{e31.8.44}). (We 
we do not use the cusp-metric structure on $\Delta\setminus\{0\}$
which comes from this identification.)

Then the $q$-fold cover $\tilde{L}^*$ of $L^*$ is given by
\begin{equation}
\tilde{L}^* =  \HH^2\times \CP^1/ q\ZZ.
\end{equation}
Since $\alpha = p/q$, the action here is trivial on $\CP^1$, and so we
can introduce coordinates
\begin{equation}\label{e32.844}
(u_1,v_1) = (e^{i\xi/q}, w_0/w_1),\;\;(u_2,v_2) = (e^{i\xi/q},
w_1/w_0), u_1=u_2,\; v_1 = v_2^{-1}.
\end{equation}
If $\omega = e^{2\pi i/q}$, then the action of $\ZZ/q\ZZ$ on
$\tilde{L}^*$ is given by the standard action of
$\Gamma_{p,q}$ on the coordinates $(u_1,v_1)$ and by 
$\Gamma_{q-p,q}$ on $(u_2,v_2)$.  Here we have written $\Gamma_{r,s}$
for the cyclic subgroup
\begin{equation}\label{e5.10.12.03}
\Gamma_{r,s} = \left\{\begin{pmatrix} e^{2\pi in/s} & 0\cr 0&
e^{2\pi i r n/s}\end{pmatrix}: n=0,1,\ldots,s-1\right\}.
\end{equation}
of $U_2$.  Thus in the orbifold $\ovL$,  there are two singularities in the fibre 
$u_1=u_2=0$.  

Denote by $\Psi: \hM \to \ovM$ the minimal resolution of
singularities. This is obtained by replacing each
$\Gamma_{p,q}$-singularity of $\ovM$ by the corresponding
Hirzebruch--Jung string (cf.\ \cite{BPVV,Fu}).  In this resolution,
$\Psi^{-1}(F_j)$ is exactly the configuration of curves constructed
in Proposition~\ref{p1.16.4.4}. Indeed, as the notation anticipates,
this resolution of singularities is isomorphic to the iterated blow-up
constructed there:

\begin{theo} Let $\chM$ be a parabolically ruled surface and let 
$M^* = \chM \backslash \cup F_j$. Let $\ovM$ by the corresponding
orbifold. Then we have the diagram
\begin{equation}\label{e10.844}
\xymatrix{&          \hM\ar_\Psi[dl] \ar^\Phi[dr]   &  \\
\ovM & & \chM \\
& M^*\ar[ul]\ar[ur]&
}
\end{equation}
where the lower arrows are the natural inclusions of $M^*$.
\label{bimero}\end{theo}

\begin{proof} 
The problem is local to the base, so we use the local models and
coordinates introduced above. In terms of these coordinates, the
identity map on $L^*$ becomes the singular map $\ovL \to \chL$
\begin{equation}\label{e33.844}
F:(u_1,v_1) \mapsto (x_1,y_1) = (u_1^q, u_1^{-p}v_1),\;
F:(u_2,v_2) \mapsto (x_2,y_2) = (u_1^q, u_1^{p}v_1),\;
\end{equation}
which is smooth away from the central fibre $u_1=u_2=0$.

We claim that this map has a very simple description in the language
of toric geometry, which makes clear the existence of the diagram
\eqref{e10.844}.

Fix the standard lattice $\ZZ^2 \subset \RR^2$. 
The toric description\footnote{Strictly, of $\CC\times \CP^1$ rather than $\Delta\times \CP^1$, but this is not important in this discussion.} 
of $\chL$ is in terms of the fan with 2-dimensional cones
$$
\sigma_1 = \{(\mu,\nu): \mu\geq 0,\nu \geq 0\},\;\;
\sigma_2 = \{(\mu,\nu): \mu\geq 0,\nu \leq 0\}.
$$
The dual cones are 
$$
\sigma_1^* = \{(m,n): m \geq 0, n\geq 0\},\;\;
\sigma_2^* = \{(m,n): m \geq 0, n\leq 0\}.
$$
The corresponding coordinate rings are just $\CC[X,Y]$ and
$\CC[X,Y^{-1}]$.  Identifying $(X,Y)$ with $(x_1,y_1)$ in the first
case and with $(x_2,y_2^{-1})$ in the second, we arrive at $\chL$,
coordinatized as in \eqref{e31.844}.
Similarly, $\ovL$ is the toric variety
corresponding to the fan with two-dimensional cones
$$
\tau_1 = \{(\mu,\nu): \mu \geq 0, p\mu + q\nu \geq  0\},\;
\tau_2 = \{(\mu,\nu): \mu \geq 0, p\mu + q\nu \leq  0\}.
$$
The dual cones are
$$
\tau_1^* = \{(m,n): qm-pn \geq 0, n\geq 0\},\;\;
\tau_2^* = \{(m,n): qm-pn \geq 0, n\leq 0\}.
$$
Using the same indeterminates as before, we have two affine varieties
with coordinate rings
$$
A_1 = \bigoplus_{qm \geq pn, n\geq 0} \CC X^mY^n,\;\;
A_2 = \bigoplus_{qm \geq pn, n\leq 0} \CC X^mY^n.
$$
Following Fulton \cite{Fu}, the first of these is identified with the 
coordinate ring of $\CC^2/\Gamma_{p,q}$ by introducing variables
$(u_1,v_1)$ with $u_1^q = X$, $Y=u_1^{-p}v$. Then $A_1$ is precisely the
$\Gamma_{p,q}$-invariant part of $\CC[u_1,v_1]$. Similarly $A_2$ is
identified with the coordinate ring of $\CC^2/\Gamma_{q-p,q}$ by
setting $X = u_2^q$, $Y = u_2^{-p}v_2^{-1}$.
Then it is clear that the identity map of $\ZZ^2$ gives rise to our
singular holomorphic map.

The map $F$ is singular because the identity map does not map the cone
$\tau_1$ into either of the cones $\sigma_1$, $\sigma_2$. The
Hirzebruch--Jung resolution of  the two singularities in $\ovL$
corresponds to a subdivision of the fan $\{\tau_1,\tau_2\}$. Recall
the continued-fraction expansions \eqref{e1.844}, \eqref{e20.844} and set
$$
v_0 = (0,1),v_1 = (1,0), v_j = (m_j,-n_j),
$$
where 
$$
\frac{n_j}{m_j} = \cfrac{1}{e_1 - \cfrac{1}{e_2-\cdots \frac{1}{e_{j-1}}}},
$$
so that $v_{k+1} = (q,-p)$. Define similarly 
$$
v_0' = (0,-1), v_1' = (1,-1), v'_j = (m'_j , n'_j - m'_j)
$$
where
$$
\frac{n'_j}{m'_j} = \cfrac{1}{e'_1 - \cfrac{1}{e'_2-\cdots
\frac{1}{e'_{j-1}}}}
$$ 
are the approximants to $(q-p)/q$. Again $v'_{l+1} = (q,-p)$. Put
$$
\rho_j = \RR_{\geq 0} v_j \oplus \RR_{\geq 0}v_{j+1},
\rho'_j = \RR_{\geq 0} v'_j \oplus \RR_{\geq 0}v'_{j+1}
$$
for every $j$. Then the fan consisting of all the $\rho_j$ and
$\rho'_j$ corresponds precisely to the minimal resolution 
$\Psi:\widehat L\to\overline L$
and the identity map of $\ZZ^2$ now induces the (restriction of the)
holomorphic map
$\Phi: \widehat L \to \check{L}$.
\end{proof}
\begin{remark*}  It is a pleasant exercise to verify from this point
of view that $\widehat L \to \check L$ 
is a multiple blow-up. Combinatorially, a
$(-1)$-curve corresponds to the occurence of adjacent cones of the
form
$$\RR_{\geq 0} u \oplus \RR_{\geq 0} (u+v),\;\;
\RR_{\geq 0} (u+v) \oplus \RR_{\geq 0} v
$$
in a fan. The blow-down operation corresponds to the replacement of
these two cones by the single cone 
$$
\RR_{\geq 0} u \oplus \RR_{\geq 0} v.
$$
In the above fan, one
has to show that there is a sequence of such deletions that always
results in a non-singular fan.  In fact, at each stage, one deletes the ray generated by the vector $v_j$ or $v'_j$ with the {\em largest} $x$-coordinate. The details are left to the interested reader.
\end{remark*}

\subsection{Scalar-flat K\"ahler orbifold metrics and holomorphic
vector fields}

\begin{theo} Let $\chM$, $\hM$ and $\ovM$ be as before, and assume
that $\chM$ is a stable parabolic ruled surface. Then $\ovM$ carries a
scalar-flat K\"ahler orbifold metric  $\ovg$; and
the algebra $\fa(\ovM)$ of  holomorphic vector fields is
$0$.
\label{theodefmbar}
\end{theo}
\begin{proof} 
Recall from \eqref{e1.154} that $\ovM$ is a quotient of $\HH^2\times \CP^1$; if the parabolic structure is stable, then $\pi_1^{orb}(\ovS)$ acts by isometries of this product, so its SFK structure descends to define an orbifold SFK structure on $\ovM$. 

The proof that $\fa(\ovM)=0$ is very close to the proof of 
\cite[Prop.~3.1]{LS}.  Suppose that
$\xi$ is a holomorphic vector field on $\ovM$.   We claim first that
$\xi$ must be vertical.  To see this, consider
the exact sequence
$$ 0 \longrightarrow \cF \longrightarrow T\overline M \stackrel  h\longrightarrow \cN
\longrightarrow 0,
$$
on $\ovM$, 
where $\cF$ is the tangent space to the fibres of $\opi$,
while $\cN$ is the normal bundle to the
fibres.  If $F=\opi^{-1}(P)$ is a smooth fibre, then $\cN|F$ is the trivial
line-bundle on $F$ and so $h(\xi)|F$ is constant on $F$ and is given by the pull-back of a vector in $T_P\ovS$. These vectors clearly patch together to give
a holomorphic vector field $V$ on the punctured Riemann surface
$\Sigma$. The same argument works near the singular fibres of $\opi$, so 
$V$ extends to define a
smooth (orbifold) vector field on $\ovS$.  But the base is an orbifold Riemann surface with negative scalar curvature, and it follows as in the smooth case that such a surface carries no non-trivial holomorphic vector fields. So $V=0$ and $\xi$ is vertical.

The vector field $\xi$ must have zeros, for it is tangent to the
smooth fibre $F$ which is a $2$-sphere.  In 
particular, it must lie in the subalgebra $\fa_0(\ovX)$ of non-parallel holomorphic vector fields, so $\fa(\ovX)= \fa_0(\ovX)$. 

Because $\ovM$ has a SFK metric $\ovg$, $\fa_0(X)$ is the complexification of the Lie algebra of non-parallel infinitesimal isometries of $\ovg$ (cf.\ \cite[Chapter 2]{Be}).  So we may assume that $\xi$ is such an infinitesimal isometry.  But then $\xi$ lifts to an infinitesimal isometry of $\HH^2\times \CP^1$, i.e. an element of the 
Lie algebra
\begin{equation}\label{e1.7.1.04}
\mathfrak{sl}_2(\RR) \times \mathfrak{su}_2
\end{equation}
which is invariant under the induced action of $\pi_1^{orb}(\ovS)$. 
This action is given by composing the representation
$$
\pi_1^{orb}(\ovS) \to  SL_2(\RR)\times SU_2/\ZZ_2
$$
with the adjoint action of $SL_2(\RR)\times SU_2/\ZZ_2$ on its Lie algebra.
Now the adjoint action of
$SU_2$ on its lie algebra is precisely the action of $SU_2/\ZZ_2 = SO_3$
on $\RR^3$, so a non-zero $\pi_1^{orb}(\ovS)$-invariant element in
\eqref{e1.7.1.04} will determine an invariant line in
$\RR^3$. An intersection of this line with the unit sphere $S^2\subset \RR^3$
defines an invariant point
 of $\CP^1\simeq S^2$, and we obtain an invariant complex line in
$\CC^2$. This contradicts the irreducibility of the representation $\rho$.
\end{proof}

\subsection{A remark on uniqueness}
\label{secunique}
Before we proceed to the proof of  Theorem~\ref{maintheo} (or
 Theorem~\ref{maintheoparab}),  it is worth mentioning that the
 SFK metrics which are produced, 
 are in fact unique in their K\"ahler class.
More generally, we have the following result which is a direct
consequence of~\cite[Theorem 1.1]{CT}. 
\begin{prop}
\label{theouniq}
Let $\widehat M$ be a complex surface as in Theorem~\ref{maintheo}.
Then, there is at most one SFK metric in each K\"ahler class.
\end{prop}
\begin{proof}
We give here an alternative proof of the proposition, based on the work
of Donaldson.
If the proposition is false, there are two distinct SFK metrics
$\omega_1$, $\omega_2$ with
the same K\"ahler class $\Omega$. The
deformation theory for SFK metrics is unobstructed \cite[Theorem
  2.8]{LS}. By density, we can find arbitrarily small SFK deformations
$\omega'_1$ and $\omega'_2$ with K\"ahler class $\Omega' \in
H^2(\widehat M,\QQ)$. We may assume that 
$\omega'_1$ and $\omega'_2$ are distinct, which is the case for
deformations small enough, and that $\Omega'$ is an integral class
after multiplication by a suitable constant factor.
By Kodaira's embedding
theorem $(\widehat M,\Omega')$ is a projective
 variety.

We remark now that $\widehat M$
cannot have any non trivial holomorphic vector
field. Let $\Xi$ be a holomorphic
vector fields on $\widehat M$.  The complex surface $\widehat M$ 
is by definition the resolution
$\pi:\widehat M \rightarrow \overline M$ of the orbifold $\overline
M$. The exceptional fibres $\pi$ have negative self-intersection hence
$\Xi$ must be tangent to them. Therefore 
the vector field $\Xi$ is projectible. By Theorem~\myref{theodefmbar} that $\overline
M$ does not admit any non trivial holomorphic vector field. Thus
$\pi_*\Xi=0$ which forces $\Xi=0$.

According to~\cite[Corollary 5]{D}, there is at most one K\"ahler metric
of constant scalar curvature on a projective manifold with no non-trivial
holomorphic vector field. Thus, $\omega'_1$ and $\omega'_2$ must be
equal. This is a contradiction, and the proposition holds.
\end{proof}

\newcommand{\bu}{\overline{u}}
\newcommand{\hg}{\widehat{g}}
\newcommand{\hu}{\widehat{u}}
\newcommand{\hnu}{\widehat{\nu}}
\newcommand{\tg}{\widetilde{g}}
\newcommand{\tz}{\widetilde{z}}
\newcommand{\ovJ}{\overline{J}}
\newcommand{\hJ}{\widehat{J}}
\newcommand{\hD}{\widehat{D}}
\newcommand{\covS}{\overline{\cS}}
\newcommand{\chS}{\widehat{\cS}}

\newcommand{\homega}{\widehat{\omega}}
\newcommand{\Om}{\Omega}
\newcommand{\om}{\omega}
\newcommand{\ci}{C^{\infty}}
\newcommand{\ve}{\varepsilon}
\newcommand{\Id}{\mathrm{Id}}

\section{Scalar-flat metrics on $\widehat{M}$}
\label{secgluing}

In this section we prove a gluing theorem that implies Theorem~\ref{maintheo} (or equivalently  Theorem~\ref{maintheoparab}). We begin with a statement of the result and an outline of the gluing argument. We then give a rapid description of the needed perturbation theory of SFK  and hermitian-ASD metrics and the necessary linear theory, before application of the implicit function theorem to prove the theorem.

\subsection{Statement of gluing theorem and outline of proof}

\begin{theo}\label{theoglue}
 Let $(\ovM, \ovom)$ be a compact SFK orbifold of complex dimension
 $2$, with all singularities isolated and cyclic. Suppose further that
 the algebra $\fa_0(\ovM)$ of non-parallel holomorphic vector fields
 on $\ovM$ is zero.  Then the minimal resolution $\widehat{M}$ of
 $\ovM$ admits a family of SFK metrics $\homega^\ve$ such that
 $\homega^\ve \to \ovom$ on any compact subset $K \subset \hM
\setminus{E}$, where $E$ is the exceptional divisor. Moreover, any
blow-up of $\hM$ carries a similar family of SFK metrics.
\end{theo}

We have just seen that the orbifold $\ovM$ associated to 
a stable parabolic ruled surface $\chM$ carries a SFK metric 
and satisfies $\fa_0(\ovM)=0$.  Our main theorems \ref{maintheo} and \ref{maintheoparab} therefore follow at once from Theorem~\ref{theoglue}.

This theorem would follow from the gluing theorems on gluing
hermitian--ASD conformal structures stated in
\cite{KS}. Unfortunately, the argument given there is not globally
consistent; we therefore give a complete proof here. While many of the
ingredients are the same, the new argument given here is perhaps a
little more direct. We shall now outline the main points;
the technical details follow in the rest of this section. 

The first step in the proof of Theorem~\ref{theoglue} is to realize
the resolution $\hM$ as a ``generalized connected sum''.  For this,
denote by $X_{p,q}$ the minimal resolution of $\CC^2/\Gamma_{p,q}$. If
$E\subset X_{p,q}$ is the exceptional divisor, then $X_{p,q}\setminus
E$ is canonically biholomorphic to $(\CC^2\setminus 0)/\Gamma_{p,q}$,
and $E$ is the Hirzebruch--Jung string constructed from the continued
fraction expansion of $p/q$ as before. The resolution of a
$\Gamma_{p,q}$-singularity at a point $0\in \ovM$ can be realized
explicitly as follows.  Choose two small positive numbers $a$ and $b$;
cut off $X_{p,q}$ at a large radius $a^{-1}$ and remove a ball
$B(0,b)$ from $\ovM$. If $z$ is an asymptotic uniformizing holomorphic
coordinate for $X_{p,q}$ and $u$ is a local uniformizing holomorphic
system near $0\in \ovM$, then we can perform the resolution by making
the identification $u = abz$.

If this model of $\hM$ is to carry an ``approximately'' SFK metric, then it is essential that $X_{p,q}$ should itself carry an asymptotically locally euclidean (ALE) SFK metric.  The existence of such metrics is guaranteed by the following:
\begin{prop}
\label{propcaldsingmet} For all relatively prime positive integers $p<q$,
$X_{p,q}$ carries an asymptotically locally euclidean SFK metric
$g$. More precisely, there exists a compact subset $K\subset X_{p,q}$
and holomorphic coordinates $z$ on the universal cover of
$X_{p,q}\backslash K$ such that
$$
|g - |\rd z|^2|  = O(|z|^{-1})
$$
for all sufficiently large $|z|$.  Moreover,
$$
|\partial^m g| = O( |z|^{-m-1})
$$
for all positive integers $m$.
\label{p1.30.1.04}\end{prop}

Such metrics were constructed in \cite{CS} for general $p$ and
$q$. Note that $X_{q-1,q}$ is the complex manifold underlying the
$A_{q-1}$ gravitational instanton of Gibbons--Hawking, 
Hitchin and Kronheimer \cite{K}.  On the other hand,
$X_{1,q}$ is the total space of $\cO(-q)\to \CP^1$, and the SFK metric
in this case is due to LeBrun
\cite{L4}. For general 
$p$ and $q$, these metrics were known to Joyce, though they do not
seem to have appeared explicitly in his published work: they are
implicit in \cite{J1} and \cite{J2}.

The proof that the metrics of \cite{CS} have the correct asymptotic properties is given in
\S\ref{s1.1.30.04}.

\vspace{12pt}

Returning to our outline of the gluing theorem, we use cut-off functions to define a sequence of approximately SFK metrics on $\hM$, by gluing the orbifold metric to the ALE metric of Proposition~\ref{propcaldsingmet}.  

Finally we use the implicit function theorem to find a genuine SFK
metric close to the approximate one. The successful completion of this
step requires in particular that a certain linear operator be surjective in a controlled way as $a$ and $b$ go to $0$.  This is why we pay a lot of attention to the linear theory in \S\ref{linth}.
The output of the implicit function theorem is only a $C^{2,\alpha}$ metric. However, any such solution will be smooth by elliptic regularity.

\vspace{12pt}

For technical reasons, we shall use the above methods to find hermitian-ASD metrics on $\hM$. By a result of Boyer \cite{Boyer}, any such conformal class must have a K\"ahler representative, and this will automatically be scalar-flat. Indeed, if 
$\omega_0$ is any representative $2$-form in the conformal class, there is a $1$-form $\beta$ defined by
\begin{equation}\label{e1.6.2.04}
\rd \omega_0 +  \beta \wedge \omega_0 = 0.
\end{equation}
Boyer proves (using the compactness of $\hM$, and the 
fact that $b_1(\hM)$ is even)
that $\beta = \rd f$, for some $f$. But then $\rd (e^{-f}\omega_0) =
0$ and we have found the desired SFK representative. 

We start our analysis now by discussing the perturbation theory of SFK metrics  and hermitian-ASD conformal structures. 

\subsection{Perturbation theory of SFK metrics}

Let $M$ be a smooth compact complex surface, with complex structure $J$. To any hermitian conformal structure $c$ on $M$ we associate the (conformal) fundamental $2$-form $\om$ as follows:
$$
c(\xi,\eta) \leftrightarrow \om(\xi,\eta) = c(\xi, J\eta)
$$
Here we think of $c$ as a positive-definite section of
$\Om^{-1/2}S^2T^*M$, where $\Om$ is the bundle of densities on $M$; accordingly, $\om$ is a weightless $(1,1)$-form, in other words, a section of the bundle
$\Om^{-1/2}\Lambda^{1,1}$.)
Fix a background conformal structure $c$ with corresponding
$(1,1)$-form $\omega$. Then the set of $J$-hermitian conformal
structures near $c$ is identified with a neighbourhood $U$ of $0$ in 
$\ci(X,\Om^{-1/2}\Lambda_0^{1,1})$; if $A\in U$, we have the new
weightless $(1,1)$-form $\omega + A$. $A$ is taken point-wise
orthogonal to $\omega$ to avoid replacing $c$ by a multiple of $c$. 
On a complex surface, $\Lambda_0^{1,1} = \Lambda^-$, so we have
parameterized the $J$-hermitian conformal structures near $c$ by a
neighbourhood of $0$ in $\ci(X,\Om^{-1/2}\Lambda^-)$.

Denote by $L_A$ the operator $\theta \mapsto (\omega +A)\wedge
\theta$, and by $\Lambda_A$ the adjoint ``trace'' map.  It is shown by
Boyer \cite{Boyer} that
$$
\Lambda_A W^+[\omega +A ] = 0 \mbox{ iff }W^+[\omega +A]=0.
$$
On the LHS we have a section of $\Lambda_{\omega+A}^+\subset
\Lambda^2$. Denote by $\cP$ the projection $\Lambda^2 \to \Lambda^+$
to the {\em fixed} subspace of $2$-forms self-dual with respect to the
fixed background $c$ and set
$$
\cF(A) =  \cP[\Lambda_A W^+[\omega + A]].
$$
(Where necessary, we shall denote the dependence of $\cF$ on the
background conformal structure $c$ by a sub- or super-script.)
The necessary facts about $\cF$ are summarized as follows:
\begin{prop}  Given a fixed $J$-hermitian conformal structure $c$,
  there is a map
$$
\cF : U \to \ci(X,\Lambda^+)
$$
where $U$ is a neighbourhood of $0 \in \ci(\Om^{-1/2}\Lambda^-)$ with
the property that $\cF^{-1}(0)$ is the set of $J$-hermitian ASD
conformal structures near $c$ on $X$. 
Furthermore, there is an expansion
\begin{equation} \label{e5.14.4.4}
\cF(A) = \cF(0) + S[A] + 
\ve_1(A, A\otimes \nabla\nabla A)+
\ve_2(A, \nabla A\otimes \nabla A)
\end{equation}
where
$$
\cF(0) =  \Lambda W^+[c],\;\;\;
S:\ci(\Om^{-1/2}\Lambda^-) \to \ci(\Lambda^+)
$$
is a conformally invariant linear elliptic operator, and the
nonlinear terms $\ve_j(A,f)$ are real-analytic in the
$0$-jet of $A$ and linear in the $0$-jet of $f$.
\label{p3.1.30.04}
\end{prop}
We note further that if $g$ is a SFK metric in the conformal class $c$
then, trivializing $\Om^{-1/2}$, the operator $S$ gets identified with the operator
\begin{equation}\label{e2.6.2.04}
S: \alpha \longmapsto \rd^+\delta\alpha + \langle\rho,\alpha\rangle\omega.
\end{equation}
where $\rho$ is the Ricci form \cite{LS}.
\subsection{Linear theory}
\label{linth}
From now on, in order to streamline the discussion,
denote by $(X_1,g_1)$ the ALE SFK space from
Proposition~\ref{p1.30.1.04} and denote by $(X_2,g_2)$ the  non-compact space
$\ovM\setminus \ovM_{sing}$. In order to save
on notation, we assume that $\ovM$ has just one singular point, giving
a conical singularity that matches the infinity of $X_1$. 
We know that there exist asymptotic (uniformizing) holomorphic
coordinates $z$, defined, say, for $|z| \geq 1/2$ and such that
$$
g_1 = |\rd z|^2 + \eta_1(z)
$$
where 
$$
|\nabla^m\eta_1|_{g_1} = O(|z|^{-m-1}).
$$
Similarly, on $X_2$, we have (uniformizing) holomorphic coordinates
$u$, say, again defined for $|u|\leq 2$ with respect to which
$$
g_2 = |\rd u|^2 + \eta_2(u),\;\;
|\eta_2(u)| = O(|u|^2),\;
|\nabla \eta_2(u)| = O(|u|),\; |\nabla^{m+2} \eta_2(u)| = O(1)
$$
for $m = 0, 1,\ldots$.

\begin{remark*} From now on we shall write, for example,
$\{|z|\leq 1\}$ as short-hand for  $X_1\setminus\{|z| >1\}$;  and similarly for
subsets of $X_2$ such as  $\{|u|\geq 1\}$.
\end{remark*}

Denote by $S_1$ and $S_2$ the linear operators from
Proposition~\ref{p3.1.30.04} determined by the metrics $g_1$ and
$g_2$.  These operators are defined over the non-compact spaces $X_1$
and $X_2$, so some care is needed in arranging for them to be
Fredholm. The needed results are now standard, having been
worked out in various forms by a number of different authors. The most
refined results can be found in Melrose's book \cite{APS}; another
useful account is \cite{LoMcO}.  This does not cover the H\"older
spaces that we shall use, however: for that one can consult Mazzeo's
paper on edge operators \cite{Ma}.

It turns out that $S_1$ and $S_2$ are Fredholm when made to operate
between suitably defined weighted H\"older (or Sobolev) spaces. In
order to define these H\"older spaces, put $$ r_1 = |z|\mbox{ for }
|z|\geq 2, r_1=1\mbox{ for }|z|\leq 1/2 $$ and $r_1\geq 1$
everywhere. Define $r_2$ similarly by continuing $$ r_2 = |u|\mbox{
for }|u|\leq 1/2 $$ to $1$ for $|u|\geq 2$, so $r_2\geq 1/2$ for
$|u|\geq 1/2$.  Define the norms $$
\|f\|_{\alpha}=\sup |f|  + \sup_{z\not=z'}\left(
 (r_1(z)+r_1(z'))^\alpha
\frac{|f(z) - f(z')|}{|z-z'|^\alpha}\right)
$$
and
$$
\|f\|_{2,\alpha}
=
\sup |f|  + \sup |r_1\nabla f| + \sup |r_1^2\nabla\nabla f|
$$
$$
\hspace{4cm}
+
\sup_{z\not=z'}\left(
 (r_1(z)+r_1(z'))^{2+\alpha}
\frac{|\nabla\nabla f(z) - \nabla \nabla f(z')|}{|z-z'|^\alpha}\right)
$$
The completion of $\ci_0(X_1)$ in $\|\cdot\|_{\alpha}$ will be
denoted by $B^\alpha(X_1)$; its completion in $\|\cdot\|_{2,\alpha}$
will be denoted by $B^{2,\alpha}$.   The weighted versions of these spaces
are
\begin{equation}\label{e1.30.1.04}
r_1^\delta B^{n,\alpha}  =
\{f: r_1^{-\delta}f \in B^{n,\alpha}(X_1)\},\;\;
\|f\|_{n,\alpha,\delta} = \|r_1^{-\delta}f\|_{n,\alpha},\;\; (n=0,2)
\end{equation}

 We define $r_2^\delta B^{n,\alpha}(X_2)$ in exactly the same way,
using $u$ and $r_2(u)$ in place of $z$ and $r_1(z)$. Similar norms
can be introduced in bundles by patching the local definitions.

The basic facts are then that
\begin{equation}\label{e3.30.1.04}
S_j : r_j^\delta B^{2,\alpha}(X_j, \Om^{-1/2}\Lambda^{-})
 \to r_j^{\delta -2}B^{0,\alpha}(X_j,\Lambda^+)
\end{equation}
is a bounded linear operator, Fredholm for all but a discrete set of
values of $\delta$. Moreover for any one of these ``good'' weights,
$S_j$ is surjective iff the formal adjoint has no null space acting on
$r_j^{-2-\delta} B^{0,\alpha}$. 

Let
\begin{equation}\label{e2.1.30.04}
S_0 = \rd^+\rd^* : \ci(\CC^2\setminus 0,\Om^{-1/2}\Lambda^-) \to
\ci(\CC^2\setminus 0, \Lambda^+)
\end{equation}
(the linearized operator at the euclidean metric).   Denote by $x$ a
standard system of euclidean coordinates on $\CC^2$.
The bad weights $\lambda$ correspond
to homogeneous solutions $\phi$, $S_0\phi =0$, 
$|\phi| = |x|^\lambda$; more precisely, they correspond to the
$\Gamma$-invariant such solutions, i.e.\ the solutions that descend to
$(\CC^2\setminus 0)/\Gamma$.

\begin{lemma} The set of bad weights is contained in $\ZZ\setminus \{-1\}$.
\label{l1.1.30.04}\end{lemma}

\begin{proof}  Suppose $\phi$ satisfies $S_0\phi=0$ and is homogeneous
of degree $\lambda\in  \CC$. By this we mean that if $\phi = \sum
\phi_je_j$, where $e_j$ is the standard parallel orthonormal basis of
$\Om^{-1/2}\Lambda^-$, then each $\phi_j$ is homogeneous of degree $\lambda$.
If $\lambda \in \CC \backslash \{-4,-5,-6,\cdots\}$, then by
\cite[Thm.~3.2.3]{Horm1}, $\phi$ has a unique extension
$\dot{\phi}$ to $\CC^2$ as a homogeneous distribution. We have
$$
S_0\dot{\phi} = f
$$
where $f$ is a distribution supported at $0$ homogeneous of degree
$\lambda - 2$, because $S_0$ is of second order.  Now any
distribution supported at $0$ is a finite linear combination of
derivatives of the $\delta$-function $\delta_0$.  Since a $k$-th order
derivative of $\delta_0$ is homogeneous of degree $-4-k$, it follows that
$$
\mbox{If }f\not=0,\; \lambda = -2-k,\mbox{ where }k=0,1,2,\ldots.
$$
On the other hand, if $f=0$, then $\dot{\phi}$ must be smooth because
$S_0$ is elliptic. In particular, the degree of homogeneity of $\phi$
must be a non-negative integer.  \end{proof}

The other essential fact about the operator $S_j$ in
\eqref{e3.30.1.04} is that any $\phi_j$ with $S_j \phi_j =0$, must have
an asymptotic expansion for $|z|\to\infty$ or $|u|\to 0$; moreover,
the leading term must behave exactly like $|z|^\lambda$ or
$|u|^\lambda$, where $\lambda$ is one of the ``bad weights'' found in
Lemma~\ref{l1.1.30.04}.

\begin{prop} If $0< \delta < 2$ and $0<\alpha < 1$ then
$$
S_1: r_1^{-\delta} B^{2,\alpha}(X_1,\Om^{-1/2}\Lambda^-)
 \to r_1^{-\delta -2}B^{0,\alpha}(X_1,\Lambda^+)
$$
has a bounded right-inverse $G_1$.  If the orbifold $\ovM$ satisfies 
$\fa_0(\ovM)=0$, then
$$
S_2: r_2^{-\delta} B^{2,\alpha}(X_2,\Om^{-1/2}\Lambda^-)
 \to r_2^{-\delta -2}B^{0,\alpha}(X_2,\Lambda^+)
$$
has a bounded right-inverse $G_2$.
\label{p6.30.1.04}\end{prop}
\begin{proof}
For $\delta$ in the given range, there are no bad weights, so $S_j$ is
Fredholm, and is surjective if and only if
$$
S_j^* \psi = 0, |\psi| = O(r_j^{\delta-2}) \Rightarrow \psi = 0.
$$
If $j=1$, then we note that
$\delta-2<0$ and so $|\psi| =
O(|z|^{-2})$.  This is enough to force vanishing of $\psi$ by
\cite[Theorem~8.4]{KS}. If $j=2$, we use that  $\delta-2 >-2$, from which
it follows that $\psi = O(1)$.  By elliptic regularity, it follows
that $\psi$ extends to a solution on the whole of the orbifold.  To
show that $\psi=0$, we follow the analysis of $S$ that was given in
\cite[Thm.~2.7]{LS}---the argument goes through without change for compact
orbifolds. The essential point is that the component of $\psi$ in the
direction of the K\"ahler form $\omega$ is a function $f$ that
satisfies Lichnerowicz's equation $$
\Delta^2 f + 2\langle \rho, \rd \rd^c f\rangle = 0.
$$
Because $\ovM$ has constant scalar curvature, $\nabla^{1,0}f \in \fa_0(\ovM)$. 
Thus $f=0$. One shows further that $f=0$ implies 
that the component of $\psi$ orthogonal to $\omega$ also vanishes.
\end{proof}

\subsection{Gluing construction}
\label{secgluingconstr}
We continue with the notation of the previous section, now picking two
small numbers $a$ and $b$.  We begin by gluing $X_1$ and $X_2$ to
produce the complex manifold $\hM = X_{a,b}$ which is the resolution of
singularities of $\ovM$.  This is easy: we just identify an annular
region $\{a^{-1}\leq |z| \leq 4a^{-1}\}$ in $X_1$ with a similar region
$\{b\leq |u|\leq 4b\}$ in $X_2$ by the holomorphic map
\begin{equation}\label{e4.30.1.04}
u  = abz.
\end{equation}

\subsubsection{Gluing metrics}
\label{glum}
The next step is to glue the metrics. 
Pick a standard cut-off function $\theta_1$, $0\leq
\theta_1 \leq 1$, with $\theta_1 =1$ for $t\leq 1$, $\theta_1 = 0$ for
$t\geq 2$. Set $\theta_2 = 1 - \theta_1$. Define new metrics $g_1^a$
on $X_1$ and $g_2^b$ on $X_2$ by 
``flattening'' $g_1$ near infinity and $g_2$ near $0$; thus $g_1^a =
g_1$ for $|z|\leq a^{-1}$ and
$$
g_1^a = |\rd z|^2  + \theta_1(a|z|)\eta_1(z)\mbox{ for }|z|\geq a^{-1}.
$$
Similarly, $g_2^b = g_2$ for $|u|\geq 4b$ and
$$
g_2^b = |\rd u|^2  + \theta_2((2b)^{-1}|u|)\eta_2(u)\mbox{ for
}|u|\leq 4b.
$$
One computes that any curvature quantity $R$ satisfies
\begin{equation}\label{er1}
|R(g_1^a)| = 
\left\{\begin{array}{ll}
|R(g_1)| & \mbox{ for } |z| \leq a^{-1} \\
O(a^{3}) & \mbox{ for } a^{-1} \leq |z| \leq 2a^{-1} \\
0 & \mbox{ for }|z| \geq   2a^{-1}.\end{array}\right.
\end{equation}
Similarly,
\begin{equation}\label{er2}
|R(g_2^b)| = 
\left\{\begin{array}{ll}
|R(g_2)| & \mbox{ for } |u| \geq 4b \\
O(1) & \mbox{ for } 2b \leq |u| \leq 4b \\
0 & \mbox{ for }|u| \leq   2b.\end{array}\right.
\end{equation}

We now define a metric on $X_{a,b}$ by matching $g_1^a$ with $a^{-2}b^{-2}g_2^b$ by
\eqref{e4.30.1.04} along the spheres $\{|u|=2a^{-1}\}$ and $\{|z| = 2b\}$:
\begin{equation}
g^{a,b} = 
\left\{\begin{array}{ll} 
 g^a_1 & \mbox{ for } |z| \leq 2a^{-1} \\
a^{-2}b^{-2}g_2^b & \mbox{ for }  |u| \geq 2b.
\end{array}\right.
\end{equation}
It is clear that $g^{a,b}$ is {\em hermitian} with respect to the
complex structure of $X_{a,b}$.

Consider now the map $\cF^{a,b}$ of Proposition~\ref{p3.1.30.04}, associated
to $g^{a,b}$ (or more accurately to the underlying conformal structure
$c^{a,b}$).  The main technical theorem can be stated as follows.
\begin{theo} There exist Banach spaces $\EE$ and $\FF$ (depending upon
$a$ and $b$) such that
\begin{enumerate}
\item $\cF:= \cF^{a,b}$ extends to a smooth map from a neighbourhood $U$ of
$0\in \EE$ to $\FF$;

\item for $x\in U$, 
$$
\cF(x) = \cF(0) + S[x] + Q(x),
$$
where $S$ is a Fredholm linear operator and $Q$ satisfies
$$
\|Q(x) - Q(y)\| \leq C_1(\|x\| + \|y\|)\|x-y\|.
$$
\item $\|\cF(0)\| \to 0$ as $a, b\to 0$.

\item If $a$ and $b$ are sufficiently small, then $S$ has a right inverse
$G$, with norm uniformly bounded by $C_2$, say.
\end{enumerate}
\label{t4.1.30.03}
\end{theo}

Once one has this result, one obtains a parameterization of
$\cF^{-1}(0)$ as a graph of a map $f$ from a small ball in $\EE_0$
into $\EE_1$, where $\EE_0$
is the finite-dimensional null-space of $S$ and $\EE_1$ is the range of
$G$.  Indeed, suppose that 
$$
\|\cF(0)\| \leq \frac{\lambda_0}{C_1C_2^2},\;\;
k_1 =\frac{\lambda_1}{C_1C_2},\;\;
k_2 =\frac{\lambda_2}{C_1C_2^2}.
$$
If we replace $x$ by $x + Gy$, where $x\in \EE_0$, then the equation
$\cF(x+Gy)=0$ becomes
$$
y = T_x(y):= - \cF(0) - Q(x+Gy)
$$
If $x$ is fixed in $\EE_0\cap\{\|x\| \leq k_1\}$ and $\|y\|\leq k_2$,
then
$$
\|T_x(y)\| \leq \|\cF(0)\| + C_1\|x+Gy\|^2
\leq \frac{\lambda_0 + \lambda_1^2 + \lambda_2^2}{C_1C_2^2}
$$
Therefore, $T_x$ maps the $k_2$-ball in $\FF$ into itself if
\begin{equation}\label{e5.30.1.04}
\lambda_0 + \lambda_1^2 + \lambda_2^2 \leq \lambda_2
\end{equation}
Furthermore,
$$
\|T_x(y) - T_x(y')\| = \|Q(x+Gy) - Q(x+Gy')\| \leq 
2(\lambda_1+\lambda_2)\|y-y'\|
$$
so that $T_x$ is a contraction mapping if
\begin{equation}\label{e6.30.1.04}
\lambda_1 + \lambda_2 < \frac{1}{2}.
\end{equation}
It is easy to find positive numbers $\lambda_0$, $\lambda_1$ and
$\lambda_2$ that simultaneously satisfy \eqref{e5.30.1.04} and
\eqref{e6.30.1.04}.  It follows from the contraction mapping theorem that
there is
a unique fixed point $y = T_x(y)$ for any given $x$ in the $k_1$-ball
of $\EE_0$, and that 
$$
\|y\| \leq \frac{1}{1 - 2(\lambda_1 +\lambda_2)}\|\cF(0)\|.
$$
The norm of the corresponding zero $x + Gy$ of $\cF$ satisfies
$$
\|x + Gy\| \leq \|x\| + \frac{C_2}{1 - 2(\lambda_1 +\lambda_2)}\|\cF(0)\|.
$$

\subsubsection{Gluing function spaces}

The functions $r_1$ and $r_2$ agree in the gluing region $b\leq |u| \leq 4b$, up to a factor of $ab$. We therefore define
\begin{equation}
w = w^{a,b} =
\left\{\begin{array}{ll} 
 r_1(z) & \mbox{ for } |z| \leq 2a^{-1} \\
a^{-1}b^{-1}r_2(u) & \mbox{ for }|u| \geq  2b.
\end{array}\right.
\end{equation}
Note that
\begin{equation}\label{e1.14.4.4}
1 \leq w \leq a^{-1}b^{-1} \mbox{ on }X^{a,b}.
\end{equation}

Set
\begin{eqnarray}
\sw\|f\|_{n,\alpha} &=&
\sup |f| + \sup|w \nabla f| + \cdots + \sup|w^n\nabla^n f| + \nonumber \\
&& \sup_{P\not= Q}(w(P) + w(Q))^{n+\alpha}
\frac{|\nabla^n f(P) - \nabla^n f(Q)|}{d(P,Q)^\alpha}.
\end{eqnarray}
Here all lengths are measured by the metric $g^{a,b}$. Denote the
completion of $\ci(X)$ in this norm by $B^{n,\alpha}(X)$. We now define
$$
\EE = w^{-\delta}B^{2,\alpha}(X,\Om^{-1/2}\Lambda^-),\;
\FF = w^{-\delta-2}B^{0,\alpha}(X,\Lambda^+).
$$
The norms in $\EE$ and $\FF$ will be denoted by subscripts $\EE$
and $\FF$ or by the notation
$$
\sw\|\cdot\|_{2,\alpha,-\delta},\;\; 
\sw\|\cdot\|_{0,\alpha,-2-\delta}.
$$
By design $\EE$ and $\FF$ are closely related to the function spaces
that were introduced in Proposition~\ref{p6.30.1.04}.  We make a couple of observations concerning $\EE$ and $\FF$.
Let  $A \in \ci(X,\Om^{1/2}\Lambda^-)$, $W\in \ci(X,\Lambda^+)$. If
$A$ and $W$ have support contained in $|z|\leq 2/a$, then clearly
\begin{equation}\label{e2.14.4.4}
\|A\|_{\EE} = \|A\|_{2,\alpha,-\delta},\;\; \|W\|_{\FF} =
\|W\|_{0,\alpha,-2-\delta}.
\end{equation}
On the other hand, if $A$ and $W$ are supported in the region $|u|\geq
2b$, then 
\begin{equation}\label{e3.14.4.4}
\|A\|_{\EE} = a^{-\delta}b^{-\delta}\|A\|_{2,\alpha,-\delta},\;\;
\|W\|_{\FF} = a^{-\delta}b^{-\delta}
\|W\|_{0,\alpha,-2-\delta}.
\end{equation}

\subsubsection*{Proof of Theorem~\ref{t4.1.30.03}} Parts (i)--(iii)
will follow from the expansion in Proposition~\ref{p3.1.30.04}.  Let
us begin by computing the $\FF$-norm of $\cF(0)$. This is zero away
from the gluing region $b \leq |u|\leq 4b$, and can be estimated there
by combining \eqref{er1} and \eqref{er2}. The result is
\begin{equation}\label{e2.22.4.4}
\|\cF(0)\|_{\FF} = O(a^{1-\delta} + a^{-\delta}b^2).
\end{equation}
Thus to guarantee (iii), we shall need to choose $0 <\delta < 1$ and
make a suitable choice of the relative sizes of $a$ and $b$, for example $a = b^2$.

In order to show that $S$ extends to a bounded linear map from $\EE$
to $\FF$, use a partition of unity to split $S(A)$ into two pieces,
the relations \eqref{e2.14.4.4} and \eqref{e3.14.4.4}, and the fact
that the operators $S_1$ and $S_2$ in Proposition~\ref{p6.30.1.04} 
are bounded; here it is important that the
same scale factor appears in comparing the norms of $A$ and $W$ in
\eqref{e3.14.4.4}.

Finally we turn to the non-linear terms in \eqref{e5.14.4.4}. These
are only defined if $\sup|A|$ is sufficiently small. But
\begin{equation}\label{e6.14.4.4}
\sup |A| \leq \sup w^\delta |A| \leq \|A\|_{\EE}
\end{equation}
by \eqref{e1.14.4.4}. Assuming, then, that $\|A\|_{\EE}$ is sufficiently small,
we have
$$
\sup|w^{2+\delta}\ve_1(A,A\otimes \nabla^2 A) \leq
\sup|\ve_1(A,w^{\delta}A\otimes w^{2+\delta}
\nabla^2 A) \leq C\|A\|^2_{\EE},
$$
again using \eqref{e1.14.4.4}. Since $\ve_1$ is linear in its second
variable, we have
$$
\ve_1(A, A\otimes\nabla^2 A) - \ve_1(B, B\otimes\nabla^2 B)
= \ve_1(A, A\otimes \nabla^2 (A - B)) +
\ve_1(A, (A-B)\otimes\nabla^2 B) +
$$
$$
\hspace{4cm} + \ve_1(A, B\otimes\nabla^2 B)
 - \ve_1(B, B\otimes\nabla^2 B)
$$
and it is straightforward to use this to prove that the H\"older
quotient in the definition of $\FF$ is controlled by $\|A\|$ and also
the quadratic estimate of $Q$ in (ii). The term $\ve_2(A,\nabla
A\otimes \nabla A)$ is estimated in the same way.

It remains to prove that $S$ has a uniformly bounded right inverse.
Note that to save on notation, we have not indicated explicitly that
$S$ depends upon $a$ and $b$. One should not lose sight of this
dependence, however, in what follows.

First note that if $a$ and $b$ are small, then the operator norms of
$S_1- S_1^a$ and $S_2- S_2^b$ are small, and so there are operators
$G_1^a$, $G_2^b$ such that
$$
S_1^a G_1^a = 1,\;\; S_2^b G_2^b = 1, \|G_1 - G_1^a\|\leq 1/2,
\|G_2 - G_2^b\|\leq 1/2.
$$
(Here the operator norms are defined by viewing $S_j$ as operators 
between the Banach spaces of Proposition~\ref{p6.30.1.04}.)
We shall now follow the approach of  
\cite[Chapter 7]{DK} to splice these right-inverses to give first an
approximate and then an exact right-inverse for $S$. 
Recall the partition of unity $\theta_1+\theta_2=1$ on $\RR$ that was
introduced  
at the beginning of \S\ref{glum}. For any small positive number
$\lambda$, which will be fixed later, define
\begin{equation}\label{e1.16.9.4}
\beta_1(z) = \theta_1\left((a|z|/4)^{\lambda}\right),\;
\beta_2(u) = \theta_2\left(2(|u|/2b)^{\lambda}\right).
\end{equation}
Define also
\begin{equation}\label{e2.16.9.4}
\gamma_1(z) = \theta_1(a|z|/2),\;\gamma_2(u) = \theta_2(|u|/2b).
\end{equation}
Then $\gamma_1 + \gamma_2=1$ on $X^{a,b}$, and
\begin{equation}\label{idemp}
\beta_1\gamma_1 = \gamma_1,\;\;\beta_2\gamma_2 = \gamma_2.
\end{equation}
Moreover, we have estimates of the form
\begin{equation}\label{e1.22.4.4}
\sup |r_1^k\nabla^k  \beta_1|=O(\lambda^k),\;\;
\sup |r_2^k\nabla^k  \beta_2| = O(\lambda^k)
\end{equation}
for each positive integer $k$, with the $O$'s uniform in $a$ and $b$.

Now form the operator on $X^{a,b}$,
\begin{equation}\label{Gdef}
G_0 = \beta_1 G_1^a \gamma_1 + \beta_2 G_2^b \gamma_2;
\end{equation}
we claim first that $\|G_0\|$, regarded as an operator $\FF \to \EE$, is
bounded independent of $a$ and $b$.  Indeed, we have 
$$
\|G_0 W \|_{\EE} \leq 
\|\beta_1 G_1^a (\gamma_1 W)\|_{\EE}
+\| \beta_2 G_2^b (\gamma_2 W)\|_{\EE}
$$
$$
\hspace{4cm}
= \|\beta_1 G_1^a (\gamma_1 W)\|_{2,\alpha,-\delta}
+
(ab)^{-\delta}\|\beta_2 G_2^b (\gamma_2 W)\|_{2,\alpha,-\delta}
$$
using the scaling formulae \eqref{e2.14.4.4} and \eqref{e3.14.4.4}. Since the operator norms of $G_1^a$ and $G_2^b$ are uniformly bounded, we obtain
$$
\|G_0 W\|_{\EE} \leq C[\|\gamma_1 W\|_{\FF} +\|\gamma_2 W\|_{\FF}]
\leq C'\|W\|_{\FF}
$$
since $\gamma_1 + \gamma_2 = 1$.

On the other hand, $G_0$ is an approximate right-inverse for $S$:
$$
SG_0 = \beta_1 S G_1^a\gamma_1
+ \beta_2 S^{a,b}G_2^a\gamma_2
+ [S^{a,b},\beta_1]G_1^a\gamma_1
+ [S^{a,b},\beta_2]G_2^a\gamma_2.
$$
On the support of $\beta_1$, $S$ is close to $S^a_1$, if $b$ is
small. Therefore,
$$
\beta_1 SG_1^a\gamma_1 = 
\beta_1 S^{a}_1G_1^a\gamma_1 +o(1) = \beta_1\gamma_1 + o(1) = \gamma_1
+o(1).
$$
(Here $o(1)$ indicates an operator whose norm tends to zero with $a$
and $b$.) Similarly,
$$
\beta_2 SG_2^b\gamma_2  =  \gamma_2 + o(1).
$$
We show now that the commutator terms are $O(\lambda)$. Clearly
$[S,\beta_1]$ is a first-order differential operator supported
where $\nabla \beta_1 \not=0$, i.e. the interval $[(4/a), (4/a)\cdot
  2^{1/\lambda}]$. For fixed $\lambda$ and small $a$, $S$ is
very close to $S_2^b$ over this interval, and we can
estimate norms as before by passing from $\EE$ and $\FF$ to the
corresponding fixed weighted H\"older spaces on $X_2$. 
The coefficients of $[S,\beta_1]$ are
therefore bounded functions on $X_2$, multiplied by $\nabla\beta_1$
and $\nabla^2\beta_1$. Hence by \eqref{e1.22.4.4}
the operator norm of $[S,\beta_1]$ is
$O(\lambda)$, as claimed. The same is true of the other commutator
term. In sum, we have found that if $\lambda$ is chosen small enough,
then as $a$ and $b$ tend to zero,
$$
SG_0 = 1 + R,\;\; \|R\| \leq 1/2
$$
Hence $G =  G_0(1+R)^{-1}$ is a controlled right-inverse of
$S$, and the proof of Theorem~\ref{t4.1.30.03}
is complete.

\vspace{10pt}
Applying the implicit function theorem as described above, we obtain
solutions in $\EE$ of our equation for all sufficiently small $a$,
$b$, subject to the constraints imposed by \eqref{e2.22.4.4}, 
and by elliptic regularity these solutions will be smooth.

Finally we claim that the corresponding SFK metric is close to the original
one. More precisely, it is close to $g_1$ on any subset of the form
$\{|z|\leq C_1\}$
of $X_1$, and it is close to $a^{-2}b^{-2}g_2$ on any subset of the
form $\{|u| \geq C_2\}$ of $X_2$.
From the implicit function theorem, our solution has fundamental 2-form
$$
\omega^{a,b} + A
$$
where $\omega^{a,b}$ is the fundamental 2-form of 
$g^{a,b}$ and
$$
\|A\|_{\EE} = O(a^{1-\delta}+ a^{-\delta}b^{2}).
$$
Choose $b^2 =a$, $\delta =1/2$, so that
$$
\sup|w^\delta A| + \sup|w^{\delta+1}\nabla A| +
\sup|w^{\delta+2}\nabla^2 A| = O(b)
$$
Then the equation
\eqref{e1.6.2.04} becomes
$$
\rd(\omega^{a,b}+A) + \beta\wedge(\omega^{a,b} + A) = 0.
$$
and so
$$
\sup w^{\delta+1}|\beta| = O(b).
$$
If we define a function $f$ with 
$\rd f = \beta$ by integrating along curves starting 
from any  base-point fixed in the region $|z| \leq 1/2$ in $X_1$,
it is not hard to see that such $f$ satisfies 
$$
\sup |f| = O(b)
$$
so that the K\"ahler form $\omega_{SFK}$ satisfies
$$
\omega_{SFK} = \omega^{a,b} + O(b)
$$
as claimed.
\subsection{Summary: the proof of Theorem~\ref{maintheoparab}}

Let us explain first how our work applies to prove that $\hM$ carries a SFK metric. By Theorem~\ref{bimero}, $\hM$ is the minimal resolution of singularities of the 
orbifold $\ovM$.  By Theorem~\ref{theodefmbar}, $\ovM$ is SFK, and $\fa_0(\ovM) = 0$. Hence by Theorem~\ref{theoglue}, $\hM$ carries SFK metrics.

Next let us note how to handle blow-ups of $\hM$.  The whole argument of this section goes through to handle ordinary blow-ups: there is a SFK metric on blow-up of the origin of 
$\CC^2$; this is the ``Burns metric'', but it also arises by taking $p=0$, $q=1$ in Proposition~\ref{propcaldsingmet}.  The only thing to check is that $\fa_0(\hM)=0$, but this follows from the fact that $\fa_0(\ovM) =0$.

\section{Asymptotics of the ALE scalar-flat K\"ahler metrics}
\label{s1.1.30.04}
This section is devoted to a proof of Proposition~\ref{p1.30.1.04}.  We
shall make extensive use of the notation of \cite{CS}. 

Recall that the metric on $X_{p,q}$ is determined by a choice of real numbers
$y_0> y_1 > y_2 > \cdots > y_k > y_{k+1}=0$, and the pairs $(m_j,n_j)$
coming from the continued fraction expansion of $p/q$. Here
$(m_{k+1},n_{k+1}) = (q,p)$ and  we define 
$(m_{k+2},n_{k+2}) = (0,1)$. 
For comparison with \cite{CS}, we note that $p$ and $q$ have
been reversed and we normalize $y_{k+1}$ to be $0$.  We shall also use
half-space coordinates $x>0$ and $y$, (rather than $(\rho,\eta)$)
so the hyperbolic metric becomes $x^{-2}(dx^2+dy^2)$.

If
$$
(a_j,b_j) = (m_j - m_{j+1}, n_j - n_{j+1})\mbox{ for }j=0,1,\ldots, k+1,
$$
then we can define
\begin{equation} \label{e1.14.11.03}
v_1 = \frac{x}{2}\sum_{j=0}^{k+1}
\frac{(a_j, b_{j})}{\sqrt{x^2 + (y - y_j)^2}}
=
\frac{x}{2\sqrt{x^2 + y^2}}(q,p-1) +
\frac{x}{2}\sum_{j=0}^k\frac{(a_j, b_{j})}{\sqrt{x^2 + (y - y_j)^2}},
\end{equation}
and
\begin{equation}\label{e2.14.11.03}
v_2 = \frac{1}{2}\sum_{j=0}^{k+1}
\frac{(y-y_j)(a_j, b_j)}{\sqrt{x^2 + (y - y_j)^2}}
=
\frac{1}{2\sqrt{x^2 + y^2}}(q,p-1) +
\frac{1}{2}\sum_{j=0}^k\frac{(y-y_j)(a_j, b_{j})}{\sqrt{x^2 + (y - y_j)^2}},
\end{equation}
giving an ALE scalar-flat K\"ahler metric
\begin{equation}\label{e3.14.11.03}
g
=\frac{ x|\langle v_1, v_2\rangle|}{x^2 + y^2}\left(
\frac{\rd x^2 + \rd y^2}{x^2}
+ \frac{\langle v_1, \rd t \rangle^2 +  \langle v_2, \rd t \rangle^2}{
\langle v_1, v_2\rangle^2}\right)
\end{equation}
where $t$ is an $T^2$-valued flat coordinate and
$\langle\cdot,\cdot\rangle$ denotes the standard symplectic form on $\RR^2$.

The asymptotic region of this metric corresponds to the point
$(x,y)=(0,0)$ in the boundary of the hyperbolic plane. We shall
analyze this metric first by introducing coordinates
\begin{equation}\label{e4.14.11.03}
R^2e^{2i\theta} = \frac{1}{y - i x}
\end{equation}
so that 
\begin{equation}\label{e5.14.11.03}
x = R^{-2}\sin 2\theta,\; y = R^{-2}\cos 2\theta,\;\;
R>0,0\leq \theta \leq \pi/2.
\end{equation}

Before attempting the computation of this metric in these coordinates,
it is worth writing the standard metric on $\CC^2$ in analogous
coordinates. Namely, if $(\tz_1,\tz_2)$ are standard complex coordinates,
let
\begin{equation}\label{e10.6.2.04}
\tz_1 = Re^{i\phi}\cos\theta,\;\;\tz_2 = Re^{i\psi}\sin\theta,\;
R>0, 0\leq \theta \leq \pi/2.
\end{equation}
Then the standard euclidean metric becomes
$$
|\rd \tz_1|^2 + |\rd \tz_2|^2 =
\rd R^2  + R^2\rd\theta^2 + R^2(\cos^2\theta \rd\phi^2
 + \sin^2\theta \rd\psi^2).
$$

Now return to the metric \eqref{e3.14.11.03}.  We shall try to
understand it for $R\gg 0$. Referring first to
\eqref{e1.14.11.03} we note that each term in the sum from $0$ to $k$
is $O(R^{-2})$, so
\begin{equation}\label{e6.14.11.03}
v_1 = \frac{\sin 2\theta}{2}(q,p-1) + O(R^{-2}).
\end{equation}
The $j$-th term in the sum in \eqref{e2.14.11.03} is
$-(a_j,b_j )+ O(R^{-2})$, so
\begin{equation} \label{e3.24.11.03}
v_2 = \frac{1}{2}(1+\cos 2\theta)(q,p-1)) + (0,1) + O(R^{-2}).
\end{equation}

Now define new angular variables
\begin{equation}\label{e7.14.11.03}
\psi = t_1/q,\; \phi =  (p/q) t_1 -t_2.
\end{equation}
The fact that the determinant is $q^{-1}$ means that $(\phi,\psi)$
really live on a $q$-fold cover.
\begin{equation}\label{e8.14.11.03}
\langle v_1,\rd t\rangle = q \sin\theta\cos\theta(\rd \psi - \rd\phi) + O(R^{-2})
\end{equation}
and
\begin{equation}\label{e9.14.11.03}
\langle v_2,\rd t\rangle = -q\sin^2\theta\,\rd\psi  - q
\cos^2\theta\,\rd\phi
 + O(R^{-2})
\end{equation}
Hence the ``angular part'' of the metric is given by
\begin{equation}
\langle v_1, \rd t \rangle^2 +  \langle v_2, \rd t \rangle^2
= q^2 \sin^2\theta \rd\psi^2  + q^2 \cos^2\theta \rd\phi^2 + O(R^{-2})
\end{equation}

On the other hand, from \eqref{e5.14.11.03},
\begin{equation}\label{e1.19.11.03}
y - i x = R^{-2}e^{-2i\theta},\;\; \rd x^2 + \rd y^2 = 4R^{-4}(\rd R^2
+ R^2\rd \theta^2).
\end{equation}
Since 
\begin{equation}\label{e2.19.11.03}
\langle v_1, v_2\rangle = \frac{q}{2}\sin 2\theta + O(R^{-2})
\end{equation}
we obtain for the metric
\begin{equation} \label{e3.19.11.03}
2q\left( \rd R^2 + R^2\rd \theta^2 + R^2[\alpha^2 
\sin^2\theta + \beta^2\cos^2\theta]\right)  + \mbox{ lower order terms
}.
\end{equation}
Here the ``lower order terms'' are just terms that are of order 2 less
than in the given expression for the metric (i.e. there are terms like
$O(R^{-2})\rd R^2$ and $O(1)\alpha^2$ etc.)

Thus the metric, pulled back to a uniformizing chart, differs from a
constant multiple of the standard flat metric by tensors of order
$R^{-2}$, as measured by the flat metric. 

\subsection{The complex structure}

The metric \eqref{e3.14.11.03} is K\"ahler with respect to the complex structure $J$,
\begin{equation}\label{e1.20.11.03}
J \rd t =
\frac{1}{\sqrt{x^2+y^2}}\left(
(xv_1 - yv_2)\frac{\rd x}{x} +  (yv_1  + x v_2)\frac{\rd y}{x}\right).
\end{equation}
It follows that the two  components of  $\rd t + iJ\rd t$ are
$(1,0)$-forms. In fact, a simple calculation shows that these
components are closed, hence holomorphic.

With respect to the coordinates introduced above, we can rewrite \eqref{e1.20.11.03}
as follows
\begin{equation}\label{e1.24.11.03}
J \rd t = - \frac{v_1}{\sin\theta \cos\theta}\rd \log R
  -  \frac{v_2}{\sin \theta \cos\theta} \rd \theta.
\end{equation}
(To check this, collect the terms in $v_1$ and $v_2$ before changing variables.)
Now with $\psi$ and $\phi$ as before, one obtains from
\eqref{e6.14.11.03} and \eqref{e3.24.11.03} that
$$
\omega_1:= \rd[ i\phi + \log R + \cos\theta] + O(1/R^2) \mbox{ and  }
\omega_2 := \rd[ i\psi + \log R + \sin\theta] + O(1/R^2);
$$
are closed holomorphic $1$-forms. Now shift to the complex coordinates
$(\tz_1,\tz_2)$, which are {\em not} $J$-holomorphic. We have
$$
\omega_j = \frac{\rd \tz_j}{\tz_j} + F_j(\tz_1,\tz_2)
$$
where the $1$-forms $F_j$ are $O(1/R^2)$.  Now set
$$
f_j(\tz_1,\tz_2) = -\int_{(\tz_1,\tz_2)}^\infty F_j
$$
where the path of integral is  $t\mapsto (t\tz_1,t\tz_2)$, $t$ from $1$ to
$\infty$. Then $\rd f_j = F_j$ and
${z}_j :=\tz_j\exp(f_j)$ are $J$-holomorphic coordinates near infinity
close to the standard ones:
$$
|z_j - \tz_j | = O(1/R)\mbox{ for }R\gg 0,
$$
since $f_j= O(1/R)$.  Re-expanding the metric in the coordinates
$(z_1,z_2)$ will add some $1/R$-terms to the expansion to
\eqref{e3.19.11.03}, but this is sufficient to complete the proof of
Proposition~\ref{p1.30.1.04}.

\end{document}